\documentclass[12pt,reqno]{amsart}
\usepackage[utf8]{inputenc}
\usepackage{tikz}
\pdfoutput=1
\usepackage{amssymb}
\usepackage{amscd}
\usepackage{amsthm}
\usepackage{amsxtra}
\usepackage[mathscr]{eucal}
\usepackage[all]{xy}
\usepackage{amsmath}
\usepackage{algorithm}
\usepackage{algpseudocode}
\usepackage{algorithmicx} 
\usepackage{algpseudocode}
\usepackage{csvsimple}
\usepackage{booktabs}
\usepackage{multicol}
\usepackage{url}
\usepackage{listings}
\usepackage{enumitem}
\usetikzlibrary{calc}

\theoremstyle{plain}
\newtheorem{theorem}{Theorem}
\newtheorem{lemma}[theorem]{Lemma}

\newtheorem{corollary}[theorem]{Corollary}

\theoremstyle{definition}
\newtheorem{definition}[theorem]{Definition}

\theoremstyle{remark}
\newtheorem{remark}[theorem]{Remark}

\theoremstyle{plain}
\newtheorem*{conjecture*}{Conjecture}
\newtheorem*{remark*}{Remark}
\newtheorem*{definition*}{Definition}

\usetikzlibrary{arrows.meta, positioning}
\usetikzlibrary{decorations.markings}
\tikzset{
  arrow/.style={
    decoration={markings,
      mark=at position 0.5 with {\arrow[scale=1.5]{stealth}}},
    postaction={decorate},
       shorten <=-2pt,
    shorten >=-2pt
  },
  vertex/.style={circle, fill=black, inner sep=2pt, outer sep=2pt}
}

\usepackage{multirow}
\usepackage{fullpage} 
\usepackage{bbm}
\usepackage{tikz-cd}
\usepackage{mathabx,amsmath,amscd, latexsym, amssymb, bbold}
\usepackage{tcolorbox}
\tcbuselibrary{breakable,xparse,skins}

\definecolor{bg}{gray}{0.95}

\usepackage{amsmath, amsthm} 
\usepackage{amsfonts}
\usepackage{amssymb}
\usepackage{fullpage} 
\usepackage{xcolor}
\usepackage{caption}
\usepackage{thm-restate}

\usepackage{mathabx,amsmath,amscd, latexsym, amssymb, bbold}
\usepackage{amsmath, amssymb, tikz}

\newcommand{\mcS}{\mathcal{S}}

\definecolor{codegreen}{rgb}{0,0.6,0}
\definecolor{codegray}{rgb}{0.5,0.5,0.5}
\definecolor{codepurple}{rgb}{0.58,0,0.82}
\definecolor{backcolour}{rgb}{0.95,0.95,0.92}

\lstdefinestyle{mystyle}{
    backgroundcolor=\color{backcolour},   
    commentstyle=\color{codegreen},
    keywordstyle=\color{magenta},
    numberstyle=\tiny\color{codegray},
    stringstyle=\color{codepurple},
    basicstyle=\footnotesize\ttfamily,
    breakatwhitespace=false,         
    breaklines=true,                 
    captionpos=b,                    
    keepspaces=true,                 
    numbers=left,                    
    numbersep=5pt,                  
    showspaces=false,                
    showstringspaces=false,
    showtabs=false,                  
    tabsize=2
}

\title{Stable Andrews--Curtis conjecture via fake surfaces and Zeeman conjecture}

\date{}

\author{Lucas Fagan}
\address{Department of Mathematics, University of California, Santa Barbara, CA 93106, USA}
\email{lucasfagan@math.ucsb.edu}

\author{Yang Qiu}
\address{Chern Institute of Mathematics and LPMC, Nankai University, Tianjin, China}
\email{yangqiu@nankai.edu.cn}

\author{Zhenghan Wang}
\address{Department of Mathematics, University of California, Santa Barbara, CA 93106, USA}
\email{zhenghwa@math.ucsb.edu}

\begin{document}

\begin{abstract}
We propose an induction scheme that aims at establishing the stable Andrews--Curtis conjecture in the affirmative. The stable Andrews--Curtis conjecture is equivalent to the conjecture that every contractible fake surface is 3-deformable to a point. We prove that every contractible fake surface of complexity less than 6 is 3-deformable to a point by induction.
\end{abstract}

\maketitle

\section{Introduction}

The stable Andrews--Curtis conjecture (sACC) is an intriguing problem closely related to the smooth 4-dimensional Poincar\'e conjecture (sPC4) \cite{andrews1965free}.  While sACC is widely believed to be false, no strong evidence has emerged. In this paper, we propose an induction scheme that aims at establishing sACC in the affirmative.  

The sACC claims that the space of balanced presentations of the trivial group is connected by several stable Andrews--Curtis moves. Any presentation $P$ of a group $\pi$ leads to a natural 2-complex $K_P$ whose fundamental group is $\pi$.  By general position, $K_P$ can be PL embedded into $\mathbb{R}^5$.  If $K_P$ is contractible, then the boundary of a closed tubular neighborhood $N(P)$ of $K_P$ in $\mathbb{R}^5$ is a PL homotopy 4-sphere.  The stable AC moves correspond to simple homotopy moves on $N(P)$, implying that if $P$ is trivial under stable AC moves, then $N(P)$ is a tubular neighborhood of a point and thus a PL 5-ball.  It follows that the boundary PL homotopy 4-sphere $\partial N(P)$ of $N(P)$ is the smooth standard 4-sphere, and thus counterexamples to sPC4 cannot be constructed in this way. 

In this paper, we study the topological version of sACC \cite{Wright1975}: that every contractible 2-complex 3-deforms to a point.  Simple homotopy by expansion and collapse is a stronger form of homotopy equivalence. The Bing house with two rooms is a well-known example of a contractible complex that cannot be collapsed to a point.

Because every 2-complex can be 3-deformed to a fake surface \cite{matveev2007algorithmic}, it suffices to prove the sACC by demonstrating that all contractible fake surfaces can be 3-deformed to a point. Our strategy to prove the sACC is to 3-deform contractible fake surfaces to a point by induction on the complexity of the fake surfaces.

We observe that contractible fake surfaces may have a weakness: one of the 2-cells in any contractible fake surface has to be embedded as conjectured in \cite{fagan2024classification}, which gives an instability in the presentations of the trivial group. This instability arises from the existence of a relator in which each generator appears at most once. This relator allows us to reduce the number of generators in the presentation by substitutions.  We do not know yet how to use this observation directly for proofs, so instead we approach the topological version by reducing the complexity of fake surfaces using embedded disks.

The Zeeman conjecture says that for any contractible 2-complex $K$, $K\times I$ can be collapsed to a point \cite{zeeman1963dunce}.  The Zeeman conjecture is believed to be false as well, but there are interesting special cases.  In the case of contractible fake surfaces that are embeddable in 3-manifolds \cite{zeeman1963dunce,gillman1983zeeman}, the Zeeman conjecture is equivalent to the 3-dimensional Poincar\'e conjecture, and therefore holds. In the case of contractible fake surfaces that are not embeddable in 3-manifolds, the Zeeman conjecture is equivalent to sACC \cite{matveev1987zeeman}. Thus by establishing sACC for low complexity, we also prove the Zeeman conjecture for these cases. 

We follow strictly the notation and terminology for fake surfaces in our earlier paper \cite{fagan2024classification}.  Fake surfaces are generic $2$-polyhedra, so they have only two kinds of singularities besides the regular manifold points: a point on an edge where three half planes join (type 1 singularity), and the center of a tetrahedron in the dual celluation of triangulations of 3-manifolds (type 2 singularity). 

We denote by $\mathcal{F}_C(s,t)$ the set of all closed contractible fake surfaces which contain $t$ vertices and whose 1-skeletons have $s$ connected components. 
Our main result is a proof of the topological version of sACC for surfaces in $\mathcal{F}_C(1,t)$, where $t\leq 5$. The case of complexity $t=1$ is the main result of \cite{ikeda1971acyclic}.  

The following is our main result, which is stated more precisely in Theorem~\ref{theorem:ac}.

\begin{theorem}
\label{thm:main}
Contractible fake surfaces of complexity less than 6 with connected 1-skeletons are all 3-deformable to a point.
\end{theorem}

\begin{corollary}
\label{cor:maxltree}
Given a contractible cellular fake surface $F=(G, \{c_i\})$ of complexity less than 6 and let $T$ be a maximal tree in $G$. If $T$ is collapsed to a base point and $P$ is the resulting presentation of the trivial group, then the stable Andrews--Curtis conjecture holds for $P$.
\end{corollary}

\begin{corollary}
Let $N(F)$ be a closed regular neighborhood of a cellular fake surface $F$ embedded in $R^5$. If the complexity of $F$ is less than 6, then $N(F)$ is the standard 5-ball $B^5$.
\end{corollary}

\begin{corollary}
The Zeeman conjecture holds for any contractible cellular fake surface $F$ of complexity less than 6.
\end{corollary}

The proof for surfaces of complexity $1$ was given by Ikeda in \cite{ikeda1971acyclic}. Our extension to complexity up to $5$ is by induction.  We give a method to reduce one fake surface to another one of lower complexity up to complexity $5$, then our result follows by induction.

\section{The Stable Andrews--Curtis Conjecture}

\begin{definition}[Stable Andrews--Curtis Moves]
Let $P = \langle x_1, \dots, x_n \mid r_1, \dots, r_n \rangle$ be a balanced group presentation. The \emph{stable Andrews--Curtis moves} are the following:
\begin{enumerate}
  \item[\textbf{(AC1)}] Replace $r_i$ with $r_ir_j$ for some $j\neq i$
  \item[\textbf{(AC2)}] Replace $r_i$ with $r_i^{-1}$
  \item[\textbf{(AC3)}] Replace $r_i$ with $gr_ig^{-1}$ where $g$ is a generator or its inverse
  \item[\textbf{(AC4)}] Add a new generator $x_{n+1}$ and trivial relator $r_{n+1}=x_{n+1}$
  \item[\textbf{(AC5)}] Remove a trivial relator and the generator, i.e., the inverse of (AC4)
\end{enumerate}
The stable Andrews--Curtis conjecture allows for two additional moves beyond the usual AC-moves (AC1)--(AC3).

Two presentations are said to be \emph{stably Andrews--Curtis equivalent} if one can be obtained from the other by a finite sequence of moves (AC1)--(AC5). 

A presentation is said to be \emph{stably Andrews--Curtis trivial} if it is stably Andrews--Curtis equivalent to the trivial presentation $\langle \,\,\, \mid \,\,\, \rangle$. 
\end{definition}

Using these moves, we can state the stable Andrews--Curtis conjecture:

\begin{conjecture*}[Stable Andrews--Curtis Conjecture]
Every balanced presentation of the trivial group is stably Andrews--Curtis trivial.
\end{conjecture*}

\bigskip
The stable Andrews--Curtis conjecture also has a topological formulation \cite{Wright1975}, which is what our approach relies on:

\begin{conjecture*}[Topological Version of sACC]
Every contractible 2-complex 3-deforms to a point.
\end{conjecture*}
Since every 2-complex 3-deforms to a fake surface, sACC is also equivalent to the conjecture that every contractible fake surface 3-deforms to a point.

\section{Embedded Disk Moves}
\label{sec:embedded-disk-moves}
In this section, we introduce the main technical tool used in this paper: two moves that 3-deform a closed fake surface with an embedded disk, depending on whether the embedded disk has a trivial or nontrivial $T$-bundle. In particular, we obtain a new surface by collapsing the union of the old surface and a 3-ball along a disk. We call these the \emph{trivial embedded disk move} and \emph{nontrivial embedded disk move}, respectively.

\medskip 

We also introduce some conventions that we use in the remainder of the paper. For a fake surface $F$, we define $\mcS_n(F)$ to be the set of type $n$ singularities of $F$. In particular, $|\mcS_2(F)|$ is the number of vertices of $F$.  

We use the notation [$x_1$,\ldots,$x_n$] to represent the attaching map of a disk, where a negative denotes inverse. 
Similarly, $(x_1,x_2,x_3)$ denotes a part of an attaching map of a disk. In Case 3, we use adjacent disks to give explicit descriptions of $T$-bundles.

\bigskip

Suppose $F_1$ is a closed fake surface with an embedded disk $D$ and there are $n$ vertices on the boundary $\partial D$ of $D$. 
\subsection{Embedded Disk with Trivial $T$-Bundle}
When the $T$-bundle over $\partial D$ is trivial, the neighborhood of $D$ in $F_1$ is as shown in Figure~\ref{fig:trivial}, where we represent the edges of 1-skeleton by green lines, vertices in $\mathcal{S}_2(F_1)$ by black points, $D$ by the blue disk. 
% We will use these notations in the following section. 
Denote by $A$ the annulus above $D$ and with one boundary component attached to $\partial D$ and the other boundary component indicated by the red circle in Figure~\ref{fig:trivial}. Set $J=D\cup A$. Then $J$ is a disk embedded in $F_1$, and $J\times I$ is a 3-ball topologically. We denote by $E$ the 3-polyhedron obtained by attaching $J\times I$ to $F_1$ such that $J\times 0=J$.

\begin{figure}
    \centering
\begin{tikzpicture}[scale=0.8]
\draw[green,thick] (0,0)..controls (1,-1) and (3,-1)..(4,0);
\draw[dashed,green,thick] (0,0)..controls (1,1) and (3,1)..(4,0);
\draw[red,thick] (0,0+2)..controls (1,-1+2) and (3,-1+2)..(4,0+2);
\draw[dashed,red,thick] (0,0+2)..controls (1,1+2) and (3,1+2)..(4,0+2);
%\draw (0,0-2)..controls (1,-1-2) and (3,-1-2)..(4,0-2);
%\draw[dashed] (0,0-2)..controls (1,1-2) and (3,1-2)..(4,0-2);

\draw[green,thick] (0,4)--(0,-4);
\draw[green,thick] (4,4)--(4,-4);
%\draw (0,4)--(-2,4);
%\draw (0,-4)--(-2,-4);
%\draw (4,4)--(6,4);
%\draw (4,-4)--(6,-4);

\filldraw (0,0) circle (2pt);
\filldraw (4,0) circle (2pt);
\node[left] at (0,0) {$A_1$};
\node[right] at (4,0) {$A_2$};

%\draw[dashed,green,thick] (2,4)--(2,-4);
%\draw[dashed] (2,4)--(2.5,4.5);
%\draw[dashed] (2,-4)--(2.5,-3.5);
%\filldraw (2,0.75) circle (2pt);
%\node[above right] at (2,0.75) {$A_3$};

\draw[dashed,green,thick] (1,4.5)--(1,-4.5);
\draw[dashed,green,thick] (3,4.5)--(3,-4.5);
%\draw[dashed] (1,4.5)--(0,5);
%\draw[dashed] (3,4.5)--(4,5);
%\draw[dashed] (1,-4.5)--(0,-5);
%\draw[dashed] (3,-4.5)--(4,-5);
%\draw[dashed,red,thick] (0,0)--(3,0.6);
\filldraw (1,0.6) circle (2pt);
\filldraw (3,0.6) circle (2pt);
\node[left] at (1,0.6) {$A_n$};
\node[right] at (3,0.6) {$A_3$};
\node at (2,1) {$......$};
\node at (2,3) {$......$};

\filldraw[opacity=0.2] (-0.1,0)--(-0.1,4)--(-2,4)--(-2,-4)--(-0.1,-4)--(-0.1,0);
\filldraw[opacity=0.2] (4.1,0)--(4.1,4)--(6,4)--(6,-4)--(4.1,-4)--(4.1,0);
\filldraw[blue,opacity=0.2] (0.1,0)..controls (1,-0.9) and (3,-0.9)..(3.9,0)..controls (3,0.9) and (1,0.9)..(0.1,0);
\filldraw[opacity=0.1] (0.9,0)--(0.9,4.5)--(0.5,5)--(0.5,-5)--(0.9,-4.5)--(0.9,0);
\filldraw[opacity=0.1] (3.1,0)--(3.1,4.5)--(3.5,5)--(3.5,-5)--(3.1,-4.5)--(3.1,0);

\filldraw[green,opacity=0.2] ;

\filldraw[red,opacity=0.1] (0.1,0+2-0.2)..controls (1,-1+2-0.1) and (3,-1+2-0.1)..(4-0.1,0+2-0.1)--(4-0.1,0+0.1-2)..controls (3,0-1+0.1-2) and (1,0-1+0.1-2)..(0.1,0.1-2)--(0.1,0+2-0.1);

\filldraw[red,opacity=0.2] (0.1,0+2-0.1)..controls (1,1+2-0.1) and (3,1+2-0.1)..(4-0.1,0+2-0.1)--(4-0.1,0+0.1-2)..controls (3,0+1+0.1-2) and (1,0+1+0.1-2)..(0.1,0.1-2)--(0.1,0+2-0.1);

\end{tikzpicture}
    \caption{Neighborhood of $D$ in $F_1$ for Trivial $T$-Bundles}
    \label{fig:trivial}
\end{figure}

\begin{figure}
    \centering
\begin{tikzpicture}[scale=0.8]
\draw[red,thick] (-1,0)..controls (1,-2) and (3,-2)..(5,0);
\draw[dashed] (-1,0)..controls (1,2) and (3,2)..(5,0);
\draw[blue,thick] (-1,0+5)..controls (1,-2+5) and (3,-2+5)..(5,0+5);
\draw (-1,0+5)..controls (1,2+5) and (3,2+5)..(5,0+5);

\draw[dashed,red,thick] (0,0)..controls (1,-1) and (3,-1)..(4,0);
\draw[dashed,blue,thick] (0,0)..controls (1,1) and (3,1)..(4,0);
\draw (0,0+5)..controls (1,-1+5) and (3,-1+5)..(4,0+5);
\draw[blue,thick] (0,0+5)..controls (1,1+5) and (3,1+5)..(4,0+5);

\draw[blue,thick] (-1,5)--(-1,0);
\draw[blue,thick] (5,5)--(5,0);
\draw[dashed,blue,thick] (0,5)--(0,0);
\draw[dashed,blue,thick] (4,5)--(4,0);

\draw[dashed,red,thick] (-1,0)--(0,0);
\draw[dashed,red,thick] (4,0)--(5,0);
\draw[blue,thick] (-1,0+5)--(0,0+5);
\draw[blue,thick] (4,0+5)--(5,0+5);

\filldraw (0,0) circle (2pt);
\filldraw (4,0) circle (2pt);
\node[above left] at (0,0) {$A_1$};
\node[above right] at (4,0) {$A_2$};

\filldraw (1,0.6) circle (2pt);
\filldraw (3,0.6) circle (2pt);
\node[above right] at (1,0.6) {$A_n$};
\node[above left] at (3,0.6) {$A_3$};
\node at (2,1) {$......$};
\node at (2,5) {$......$};

\draw[dashed] (1,0.6)--(1,0.6+5);
\draw[dashed] (3,0.6)--(3,0.6+5);
\draw[dashed] (1,0.6)--(0.5,1.15);
\draw[dashed] (0.5,1.15)--(0.5,1.15+5);
\draw[dashed] (3,0.6)--(3.5,1.15);
\draw[dashed] (3.5,1.15)--(3.5,1.15+5);
\draw[dashed] (1,0.6+5)--(0.5,1.15+5);
\draw[dashed] (3,0.6+5)--(3.5,1.15+5);

\begin{scope}[shift={(8,0)}]
\draw (-1,0)..controls (1,-2) and (3,-2)..(5,0);
\draw[dashed,red,thick] (-1,0)..controls (1,2) and (3,2)..(5,0);
\draw (-1,0+5)..controls (1,-2+5) and (3,-2+5)..(5,0+5);
\draw[blue,thick] (-1,0+5)..controls (1,2+5) and (3,2+5)..(5,0+5);

\draw[dashed] (0,0)..controls (1,-1) and (3,-1)..(4,0);
\draw[dashed,red,thick] (0,0)..controls (1,1) and (3,1)..(4,0);
\draw (0,0+5)..controls (1,-1+5) and (3,-1+5)..(4,0+5);
\draw[blue,thick] (0,0+5)..controls (1,1+5) and (3,1+5)..(4,0+5);

\draw[blue,thick] (-1,5)--(-1,0);
\draw[blue,thick] (5,5)--(5,0);
\draw[dashed,blue,thick] (0,5)--(0,0);
\draw[dashed,blue,thick] (4,5)--(4,0);

\draw[dashed,red,thick] (-1,0)--(0,0);
\draw[dashed,red,thick] (4,0)--(5,0);
\draw[blue,thick] (-1,0+5)--(0,0+5);
\draw[blue,thick] (4,0+5)--(5,0+5);

\filldraw (0,0) circle (2pt);
\filldraw (4,0) circle (2pt);
\node[above left] at (0,0) {$A_1$};
\node[above right] at (4,0) {$A_2$};

\filldraw (1,0.6) circle (2pt);
\filldraw (3,0.6) circle (2pt);
\node[above right] at (1,0.6) {$A_n$};
\node[above left] at (3,0.6) {$A_3$};
\node at (2,1) {$......$};
\node at (2,5) {$......$};

\draw[dashed,blue,thick] (1,0.6)--(1,0.6+5);
\draw[dashed,blue,thick] (3,0.6)--(3,0.6+5);
\draw[dashed,red,thick] (1,0.6)--(0.5,1.15);
\draw[dashed,blue,thick] (0.5,1.15)--(0.5,1.15+5);
\draw[dashed,red,thick] (3,0.6)--(3.5,1.15);
\draw[dashed,blue,thick] (3.5,1.15)--(3.5,1.15+5);
\draw[dashed,blue,thick] (1,0.6+5)--(0.5,1.15+5);
\draw[dashed,blue,thick] (3,0.6+5)--(3.5,1.15+5);
\end{scope}
\end{tikzpicture}
    \caption{$J\times I$ for Trivial $T$-Bundles}
    \label{fig:trivialJtimesI}
\end{figure}

\begin{figure}
    \centering
\begin{tikzpicture}[scale=0.8]
\draw (-1,0)..controls (1,-2) and (3,-2)..(5,0);
\draw[dashed] (-1,0)..controls (1,2) and (3,2)..(5,0);
\draw (-1,0+5)..controls (1,-2+5) and (3,-2+5)..(5,0+5);
\draw (-1,0+5)..controls (1,2+5) and (3,2+5)..(5,0+5);

\draw[dashed,green,thick] (0,0)..controls (1,1) and (3,1)..(4,0);
\draw[dashed] (0,0)..controls (1,-1) and (3,-1)..(4,0);
\draw[green,thick] (0,0+5)..controls (1,1+5) and (3,1+5)..(4,0+5);

\draw[green,thick] (-1,5)--(-1,0);
\draw[green,thick] (5,5)--(5,0);

\draw[dashed] (-1,0)--(0,0);
\draw[dashed] (4,0)--(5,0);
\draw[green,thick] (-1,0+5)--(0,0+5);
\draw[green,thick] (4,0+5)--(5,0+5);

\draw[green,thick] (-1,0)--(-1,-2);
\draw[green,thick] (5,0)--(5,-2);
\draw[dashed,green,thick] (0,0)--(0,-2);
\draw[dashed,green,thick] (4,0)--(4,-2);

\filldraw (1,0.6) circle (2pt);
\filldraw (3,0.6) circle (2pt);
\node[below right] at (1,0.6) {$A_n$};
\node[below left] at (3,0.6) {$A_3$};
\node at (2,1) {$......$};
\node at (2,5) {$......$};

\draw[dashed,green,thick] (1,0.6)--(1,0.6+5);
\draw[dashed,green,thick] (3,0.6)--(3,0.6+5);
\draw[dashed] (1,0.6)--(0.5,1.15);
\draw[dashed,green,thick] (0.5,1.15)--(0.5,1.15+5);
\draw[dashed] (3,0.6)--(3.5,1.15);
\draw[dashed,green,thick] (3.5,1.15)--(3.5,1.15+5);
\draw[green,thick] (1,0.6+5)--(0.5,1.15+5);
\draw[green,thick] (3,0.6+5)--(3.5,1.15+5);
\draw[dashed,green,thick] (1,0.6)--(1,0.6-0.5);
\draw[dashed,green,thick] (0.5,1.15)--(0.5,1.15-0.5);
\draw[dashed,green,thick] (3,0.6)--(3,0.6-0.5);
\draw[dashed,green,thick] (3.5,1.15)--(3.5,1.15-0.5);

\filldraw (1,0.6+5) circle (2pt);
\filldraw (3,0.6+5) circle (2pt);
\node[above right] at (1,0.6+5) {$A_n^{\prime}$};
\node[above left] at (3,0.6+5) {$A_3^{\prime}$};

\node[right] at (0,0) {$A_1$};
\node[left] at (4,0) {$A_2$};

\node at (2,-3) {$n\geq3$};

\begin{scope}[shift={(8,0)}]
\draw (-1,0)..controls (1,-2) and (3,-2)..(5,0);
\draw[dashed] (-1,0)..controls (1,2) and (3,2)..(5,0);
\draw (-1,0+5)..controls (1,-2+5) and (3,-2+5)..(5,0+5);
\draw (-1,0+5)..controls (1,2+5) and (3,2+5)..(5,0+5);

\draw[dashed,green,thick] (0,0)..controls (1,1) and (3,1)..(4,0);
\draw[dashed] (0,0)..controls (1,-1) and (3,-1)..(4,0);
\draw[green,thick] (0,0+5)..controls (1,1+5) and (3,1+5)..(4,0+5);

\draw[green,thick] (-1,5)--(-1,0);
\draw[green,thick] (5,5)--(5,0);

\draw[dashed] (-1,0)--(0,0);
\draw[dashed] (4,0)--(5,0);
\draw[green,thick] (-1,0+5)--(0,0+5);
\draw[green,thick] (4,0+5)--(5,0+5);

\draw[green,thick] (-1,0)--(-1,-2);
\draw[green,thick] (5,0)--(5,-2);
\draw[dashed,green,thick] (0,0)--(0,-2);
\draw[dashed,green,thick] (4,0)--(4,-2);

\node[right] at (0,0) {$A_1$};
\node[left] at (4,0) {$A_2$};

\node at (2,-3) {$n=2$};
\end{scope}

\end{tikzpicture}
    \caption{$F_2$ for Trivial $T$-Bundles}
    \label{fig:F_2trivial}
\end{figure}

\begin{figure}
    \centering
\begin{tikzpicture}[scale=0.8]
\draw (-1,0)..controls (1,-2) and (3,-2)..(5,0);
\draw[dashed] (-1,0)..controls (1,2) and (3,2)..(5,0);
\draw (-1,0+5)..controls (1,-2+5) and (3,-2+5)..(5,0+5);
\draw (-1,0+5)..controls (1,2+5) and (3,2+5)..(5,0+5);

\draw[dashed] (0,0)..controls (1,-1) and (3,-1)..(4,0);
\draw[dashed] (0,0)..controls (1,1) and (3,1)..(4,0);

\draw (-1,5)--(-1,0);
\draw (5,5)--(5,0);

\draw[dashed,red,thick] (0,0)--(-1,0);
\draw[dashed,green,thick] (0,0)--(0,-1.5);
\draw[green,thick] (-1,0)--(-1,-1.5);

\draw[->] (-0.3,0)--(-0.3,-1.5);
\draw[->] (-0.7,0)--(-0.7,-1.5);
\node[below] at (-0.5,-1.5) {\text{collapsing}};

\draw[dashed] (4,0)--(4,-1.5);

\filldraw[opacity=0.2] (-0.1,-0.1)--(-0.9,-0.1)--(-0.9,-1.5)--(-0.1,-1.5)--(-0.1,-0.1);
\end{tikzpicture}
    \caption{$F_2$ with Free 1-Face for $n=1$}
    \label{fig:F_2free}
\end{figure}

As shown in Figure~\ref{fig:trivialJtimesI}, $J\times I$ consists of $n$ 3-subpolyhedrons whose 1-faces are indicated by red and blue lines. They are all 3-balls topologically. Now we collapse these 3-balls along certain 2-faces to get a new closed fake surface $F_2$. First, we collapse the 3-ball indicated by the left diagram in Figure~\ref{fig:trivialJtimesI} along the free 2-face enclosed by red lines. Next we collapse the other 3-balls indicated by the right diagram in Figure~\ref{fig:trivialJtimesI} by pushing upwards the 2-faces enclosed by red lines. We call this sequence of 3-deformations the \emph{trivial embedded disk move}.

When $n\geq2$, we get a new closed fake surface locally as shown in Figure~\ref{fig:F_2trivial}. The original vertices $A_1$ and $A_2$ are removed and new vertices $A_3^{\prime},\ldots,A_n^{\prime}$ are introduced, so the new surface $F_2$ contains $|\mathcal{S}_2(F_1)|+n-4$ vertices. 

In particular, for $n=2$, the number of vertices decreases by 2: $|\mcS_2(F_2)|=|\mcS_2(F_1)|-2$. Green edges are separated into two disjoint parts locally as shown in Figure~\ref{fig:F_2trivial}. Thus $F_2$ may have one more connected component in its 1-skeleton than $F_1$. For $n=3$, the number of vertices decreases by 1. Green edges are still connected together, so $F_2$ has the same number of connected components of 1-skeleton as $F_1$. In this case, the move is simply the $T^{-1}$ move from \cite{matveev2007algorithmic}.

When $n=1$, the collapses above introduce a free 1-face as indicated by the red line in Figure~\ref{fig:F_2free}. Collapsing the disk containing this free 1-face may increase the number of connected components of the 1-skeleton. If the complexity is greater than 1, the embedded disk $D$ must contain a vertex other than $A_1$. Thus $|\mathcal{S}_2(F_2)|\leq\text{max}\{|\mathcal{S}_2(F_1)|-2,0\}$.  

The effect of this move is shown in Figure~\ref{fig:trivial-embedded-disk-move}. The maximal case in which the complexity is reduced, i.e., when $n=3$, is shown in Figure~\ref{fig:trivial-embedded-disk-move-3-verts}. As mentioned above, in this case the move is equivalent to the $T^{-1}$ move from \cite{matveev2007algorithmic}. 

\begin{figure}
    \centering
\begin{tikzpicture}[scale=1.5, thick]
  % First image (left)
  \begin{scope}
  % Define coordinates
  \coordinate (A) at (0,0);
  \coordinate (B) at (4,0);

  % Fill the region between the two curved edges in grey
  \fill[green!30] (A) to[bend left=40] (B) to[bend left=40] (A);

  % Define coordinates along the top arc for the 4 vertices
  \path[bend left=40] (A) to coordinate[pos=0.2] (V1)
                           coordinate[pos=0.4] (V2)
                           coordinate[pos=0.6] (V3)
                           coordinate[pos=0.8] (V4) (B);

  % Draw curved arcs
  \draw[bend left=40] (A) to (B);
  \draw[red, very thick, bend right=40] (A) to (B);

  % Draw edges from each vertex (upward and downward)
  \draw[red, very thick] (A) -- ++(0,0.8);
  \draw (A) -- ++(0,-0.8);
  \draw[red, very thick] (B) -- ++(0,0.8);
  \draw (B) -- ++(0,-0.8);
  \draw (V1) -- ++(0,0.8);
  \draw (V1) -- ++(0,-0.8);
  \draw (V2) -- ++(0,0.8);
  \draw (V2) -- ++(0,-0.8);
  \draw (V3) -- ++(0,0.8);
  \draw (V3) -- ++(0,-0.8);
  \draw (V4) -- ++(0,0.8);
  \draw (V4) -- ++(0,-0.8);

    % Draw vertices
  \fill (A) circle (2pt) node[below] {};
  \fill (B) circle (2pt) node[below] {};

  % Draw the vertex circles (after edges so they appear on top)
  \fill (V1) circle (2pt);
  \fill (V2) circle (2pt);
  \fill (V3) circle (2pt);
  \fill (V4) circle (2pt);

  % Optionally, add labels to the arcs
  \node at (2,0.9) {};
  \node at (2,-0.9) {};
  \end{scope}

  % add an arrow from the first image to the second
    \draw[->, thick] (4.5,0) -- (5.5,0);

  % Second image (right) - copy of the first
  \begin{scope}[xshift=6cm]
  % Define coordinates
  \coordinate (A) at (0,0);
  \coordinate (B) at (4,0);

  % Draw vertices
%   \fill (A) circle (2pt) node[below] {};
%   \fill (B) circle (2pt) node[below] {};

  % Define bottom points below A and B
  \coordinate (BottomA) at (0,-0.8);
  \coordinate (BottomB) at (4,-0.8);

  % Fill the region: top edge (A to B), down from B, across bottom, up to A
%   \fill[gray!30] (A) to[bend left=40] (B) -- (BottomB) -- (BottomA) -- cycle;

  % Define coordinates along the top arc for the 4 vertices
  \path[bend left=40] (A) to coordinate[pos=0.2] (V1)
                           coordinate[pos=0.4] (V2)
                           coordinate[pos=0.6] (V3)
                           coordinate[pos=0.8] (V4) (B);

  % Define top vertices
  \coordinate (T1) at (V1 |- 0,1.5);
  \coordinate (T2) at (V2 |- 0,1.5);
  \coordinate (T3) at (V3 |- 0,1.5);
  \coordinate (T4) at (V4 |- 0,1.5);

  % Fill region between V1 and V4, from top curve to horizontal line
  \fill[green!30] (V1) -- (T1) -- (T2) -- (T3) -- (T4) -- (V4) 
    to[bend right=21] (V1) -- cycle;

  % Draw curved arcs
  \draw[bend left=40] (A) to (B);
%   \draw[bend right=40] (A) to (B);

  % Draw edges from each vertex (upward and downward)
  \draw (A) -- ++(0,-0.8);
  \draw (B) -- ++(0,-0.8);
  \draw (V1) -- (V1 |- 0,1.5);
  \draw (V1) -- ++(0,-0.8);
  \draw (V2) -- (V2 |- 0,1.5);
  \draw (V2) -- ++(0,-0.8);
  \draw (V3) -- (V3 |- 0,1.5);
  \draw (V3) -- ++(0,-0.8);
  \draw (V4) -- (V4 |- 0,1.5);
  \draw (V4) -- ++(0,-0.8);

  % Define top vertices
  \coordinate (T1) at (V1 |- 0,1.5);
  \coordinate (T2) at (V2 |- 0,1.5);
  \coordinate (T3) at (V3 |- 0,1.5);
  \coordinate (T4) at (V4 |- 0,1.5);

  % Connect top vertices with horizontal lines
  \draw[red, very thick] (T1) -- (0,1.5);
  \draw[red, very thick] (T4) -- (4,1.5);
  \draw[red, very thick] (T1) -- (T2);
  \draw[red, very thick] (T2) -- (T3);
  \draw[red, very thick] (T3) -- (T4);

  % Draw the vertex circles (after edges so they appear on top)
  \fill (V1) circle (2pt);
  \fill (V2) circle (2pt);
  \fill (V3) circle (2pt);
  \fill (V4) circle (2pt);

    \draw (T1) -- ++(0,0.5);
    \draw (T2) -- ++(0,0.5);
    \draw (T3) -- ++(0,0.5);
    \draw (T4) -- ++(0,0.5);

  % Draw top vertex circles
  \fill[red] (T1) circle (2pt);
  \fill[red] (T2) circle (2pt);
  \fill[red] (T3) circle (2pt);
  \fill[red] (T4) circle (2pt);

  % Optionally, add labels to the arcs
  \node at (2,0.9) {};
  \node at (2,-0.9) {};
  \end{scope}
\end{tikzpicture}
    \caption{The trivial embedded disk move. Embedded disks are colored, with green indicating they have trivial $T$-bundles. After applying the move, two vertices are removed, and $n-2$ new vertices are created.}
    \label{fig:trivial-embedded-disk-move}
\end{figure}

\begin{figure}
    \centering
    \begin{tikzpicture}[scale=1.5, thick]
  % First image (left)
  \begin{scope}
  % Define coordinates
  \coordinate (A) at (0,0);
  \coordinate (B) at (4,0);

  % Fill the region between the two curved edges in grey
  \fill[green!30] (A) to[bend left=40] (B) to[bend left=40] (A);

  % Define coordinates along the top arc for the 4 vertices
  \path[bend left=40] (A) to coordinate[pos=0.5] (V1)
                        (B);

  % Draw curved arcs
  \draw[bend left=40] (A) to (B);
  \draw[red, very thick, bend right=40] (A) to (B);

  % Draw edges from each vertex (upward and downward)
  \draw[red, very thick] (A) -- ++(0,0.8);
  \draw (A) -- ++(0,-0.8);
  \draw[red, very thick] (B) -- ++(0,0.8);
  \draw (B) -- ++(0,-0.8);
  \draw (V1) -- ++(0,0.8);
  \draw (V1) -- ++(0,-0.8);

    % Draw vertices
  \fill (A) circle (2pt) node[below] {};
  \fill (B) circle (2pt) node[below] {};

  % Draw the vertex circles (after edges so they appear on top)
  \fill (V1) circle (2pt);

  % Optionally, add labels to the arcs
  \end{scope}

  % add an arrow from the first image to the second
    \draw[->, thick] (4.5,0) -- (5.5,0);
    \node at (5,0.3) {$T^{-1}$};

  % Second image (right) - copy of the first
  \begin{scope}[xshift=6cm]
  % Define coordinates
  \coordinate (A) at (0,0);
  \coordinate (B) at (4,0);

  % Draw vertices
%   \fill (A) circle (2pt) node[below] {};
%   \fill (B) circle (2pt) node[below] {};

  % Define bottom points below A and B
  \coordinate (BottomA) at (0,-0.8);
  \coordinate (BottomB) at (4,-0.8);

  % Fill the region: top edge (A to B), down from B, across bottom, up to A
%   \fill[gray!30] (A) to[bend left=40] (B) -- (BottomB) -- (BottomA) -- cycle;

  % Define coordinates along the top arc for the 4 vertices
  \path[bend left=40] (A) to coordinate[pos=0.5] (V1) (B);

  % Define top vertices
  \coordinate (T1) at (V1 |- 0,1.5);

  % Fill region between V1 and V4, from top curve to horizontal line
%   \fill[green!30] (V1) -- (T1) -- (T2) -- (T3) -- (T4) -- (V4) 
    % to[bend right=21] (V1) -- cycle;

  % Draw curved arcs
  \draw[bend left=40] (A) to (B);
%   \draw[bend right=40] (A) to (B);

  % Draw edges from each vertex (upward and downward)
  \draw (A) -- ++(0,-0.8);
  \draw (B) -- ++(0,-0.8);
  \draw (V1) -- (V1 |- 0,1.5);
  \draw (V1) -- ++(0,-0.8);

  % Define top vertices
  \coordinate (T1) at (V1 |- 0,1.5);

  % Connect top vertices with horizontal lines
  \draw[red, very thick] (T1) -- (0,1.5);
  \draw[red, very thick] (T1) -- (4,1.5);

  % Draw the vertex circles (after edges so they appear on top)
  \fill (V1) circle (2pt);

    \draw (T1) -- ++(0,0.5);

  % Draw top vertex circles
  \fill[red] (T1) circle (2pt);

  \end{scope}
\end{tikzpicture}

    \caption{The trivial embedded disk move in the case $n=3$. This is the maximal case in which the complexity is reduced. In this case, the move is equivalent to the $T^{-1}$ move from \cite{matveev2007algorithmic}.}
    \label{fig:trivial-embedded-disk-move-3-verts}
\end{figure}

\subsection{Embedded Disk with Nontrivial $T$-Bundle}

When the $T$-bundle over $\partial D$ is nontrivial, the neighborhood of $D$ in $F_1$ is as shown in Figure~\ref{fig:nontrivial}. Denote by $B$ the upper band of the M\"obius band, where one long edge is attached to $\partial D$ and the other long edge is indicated by the red line in the bottom diagram in Figure~\ref{fig:Jnontrivial}. Set $J=D\cup B$. Then $J$ is an embedded disk in $F_1$ as shown in Figure~\ref{fig:Jnontrivial}. Attach 3-ball $J\times I$ to $F_1$ along $J\times 0=J$ to get $E$.

\begin{figure}
    \centering
\begin{tikzpicture}[scale=0.8]
\draw[green,thick] (0,0)..controls (1,-1) and (3,-1)..(4,0);
\draw[dashed,green,thick] (0,0)..controls (1,1) and (3,1)..(4,0);
\draw[red,thick] (0,0+2)..controls (1,-1+2) and (3,-1+2)..(4,0+2);
\draw[dashed,red,thick](4,0+2)..controls (3,2.6)..(2,2.7);
\draw[dashed](4,0+2-4)..controls (3,2.6-4)..(2,2.7-4);
%\draw[dashed] (0,0+2)..controls (1,1+2) and (3,1+2)..(4,0+2);
\draw (0,0-2)..controls (1,-1-2) and (3,-1-2)..(4,0-2);
%\draw[dashed] (0,0-2)..controls (1,1-2) and (3,1-2)..(4,0-2);

\draw[green,thick] (0,4)--(0,-4);
\draw[green,thick] (4,4)--(4,-4);
\draw (0,4)--(-2,4);
\draw (0,-4)--(-2,-4);
\draw (4,4)--(6,4);
\draw (4,-4)--(6,-4);

\filldraw (0,0) circle (2pt);
\filldraw (4,0) circle (2pt);
\node[left] at (0,0) {$A_1$};
\node[right] at (4,0) {$A_2$};

%\draw[dashed,green,thick] (1,5)--(1,-5);
\draw[dashed,green,thick] (3,5)--(3,-5);
%\draw[dashed] (1,5)--(-1.5,6);
\draw[dashed] (3,5)--(5,6);
%\draw[dashed] (1,-5)--(-1.5,-6);
\draw[dashed] (3,-5)--(5,-6);

\draw[dashed,red,thick] (2,2.7)..controls(1,1)..(0,-2);
\draw[dashed] (2,2.7-4)--(1,0.5);
\draw[dashed] (0.8,0.7)--(0,2);
%\draw[dashed] (2,2.7-4)--(0,2);

%\filldraw (1,0.6) circle (2pt);
\filldraw (3,0.6) circle (2pt);
%\node[left] at (1,0.6) {$A_n$};
\node[right] at (3,0.6) {$A_n$};

\node at (3.5,1) {\rotatebox{150}{$......$}};
\node at (4,4.6) {\rotatebox{150}{$......$}};

\end{tikzpicture}     
    \caption{Neighborhood of $D$ in $F_1$ for Nontrivial $T$-Bundles}    \label{fig:nontrivial}
\end{figure}

\begin{figure}
    \centering
\begin{tikzpicture}[scale=0.8]
\draw[green,thick] (0,0)..controls (1,-1) and (3,-1)..(4,0);
\draw[dashed,green,thick] (0,0)..controls (1,1) and (3,1)..(4,0);
\draw[red,thick] (0,0+2)..controls (1,-1+2) and (3,-1+2)..(4,0+2);
\draw[red,thick] (4,0+2)..controls (3,2.6)..(2,2.7);

\draw[green,thick] (0,2)--(0,-2);
\draw[green,thick] (4,0)--(4,2);

\filldraw (0,0) circle (2pt);
\filldraw (4,0) circle (2pt);
\node[left] at (0,0) {$A_1$};
\node[right] at (4,0) {$A_2$};

\draw[red,thick] (2,2.7)--(1.15,1.5);
\draw[dashed,red,thick] (1,1.3)--(0.5,-0.2);
\draw[red,thick] (0.4,-0.4)--(0,-2);

\filldraw (3,0.6) circle (2pt);

\node[right] at (3,0.6) {$A_n$};

\node at (3.5,1) {\rotatebox{150}{$......$}};

\draw[dashed,green,thick] (3,0.6)--(3,2.5);

\filldraw[blue,opacity=0.2] (0.1,0)..controls (1,-0.9) and (3,-0.9)..(3.9,0)..controls (3,0.9) and (1,0.9)..(0.1,0);

\filldraw[red,opacity=0.2] (0.1,0)..controls (1,-0.9) and (3,-0.9)..(3.9,0)--(3.9,1.9)..controls (3,-0.9+2) and (1,-0.9+2)..(0.1,1.9)--(0.1,0);

\filldraw[red,opacity=0.2] (3.9,1.9)..controls (3,2.5)..(2,2.6)--(1.2,1.4)--(0.5,-0.2)--(0.1,-1.5)--(0.1,0.1)..controls (1,1.1) and (3,1.1)..(3.9,0.1)--(3.9,1.9);

\begin{scope}[shift={(8,0)}]

\draw[green, thick] (2,0) circle (1);
\draw[red,thick] (4,0) arc(0:150:2);
\draw[red,thick] (4,0) arc(0:-150:2);
\draw[green,thick] (0.25,1)--(1,0);
\draw[green,thick] (1,0)--(0.25,-1);

\node[left] at (1,0){$A_1$};
\node[above right] at (3,0){$A_2$};
\node[above right] at (2,1){$A_n$};

\filldraw (1,0) circle(2pt);
\filldraw (3,0) circle(2pt);
\filldraw (2,1) circle(2pt);

\draw[green,thick] (3,0)--(4,0);
\draw[green,thick] (2,1)--(2,2);

\node at (3,1) {\rotatebox{140}{$......$}};

\filldraw[blue,opacity=0.2] (2,0) circle (0.9);
\filldraw[red,opacity=0.2] (1,0)--(0.25,-1) arc(-150:150:2)--(1,0) arc(180:-180:1);
\end{scope}

\begin{scope}[shift={(3,-4)}]
\draw[green,thick] (0,0)--(6,0);
\draw[green,thick] (0,-1)--(0,1);
\draw[green,thick] (6,-1)--(6,1);
\draw[red,thick] (0,1)--(6,1);
\draw (0,-1)--(6,-1);
\draw[green,thick] (1,-1)--(1,1);
\draw[green,thick] (3,-1)--(3,1);

\filldraw (0,0) circle (2pt);
\filldraw (1,0) circle (2pt);
\filldraw (3,0) circle (2pt);
\filldraw (6,0) circle (2pt);

\node[left] at (0,0) {$A_1$};
\node[right] at (6,0) {$A_1$};
\node[below right] at (1,0) {$A_2$};
\node[below right] at (3,0) {$A_n$};

\node[above] at (2,0) {......};

\draw[-{Stealth[length=6pt]}] (0,0.6) -- ++(90:0.01);
% \node at (0,0.5) {\rotatebox{90}{$>$}};
\draw[-{Stealth[length=6pt]}] (0,-0.3) -- ++(90:0.01);
\draw[-{Stealth[length=6pt]}] (0,-0.55) -- ++(90:0.01);

% \node at (0,-0.5) {\rotatebox{90}{$>>$}};
% \node at (6,-0.5) {\rotatebox{270}{$>$}};
\draw[-{Stealth[length=6pt]}] (6,-0.7) -- ++(270:0.01);
\draw[-{Stealth[length=6pt]}] (6,0.2) -- ++(270:0.01);
\draw[-{Stealth[length=6pt]}] (6,0.45) -- ++(270:0.01);

% \node at (6,0.5) {\rotatebox{270}{$>>$}};

\filldraw[red,opacity=0.2] (0,0)--(6,0)--(6,1)--(0,1)--(0,0);

\end{scope}

\end{tikzpicture}
    \caption{$J$ for Nontrivial $T$-Bundles}    \label{fig:Jnontrivial}
\end{figure}

\begin{figure}
    \centering
\begin{tikzpicture}[scale=0.8]
\draw[dashed,blue,thick] (0,0)..controls (1,1) and (2,1)..(3,0);
\draw[dashed,red,thick] (0,0)..controls (1,-1) and (2,-1)..(3,0);

\draw[dashed] (-2,0)..controls (1,2) and (2,2)..(4,0);
\draw[red,thick] (-1,-1)..controls (1,-2) and (2,-2)..(4,0);

\draw (-2,0)--(-1.1,0);
\draw (-0.9,0)--(0,0);
\draw[dashed,red,thick] (-1,-1)--(0,0);

\draw[blue,thick] (0,0+4)..controls (1,1+4) and (2,1+4)..(3,0+4);
\draw (0,0+4)..controls (1,-1+4) and (2,-1+4)..(3,0+4);

\draw (-2,0+4)..controls (1,2+4) and (2,2+4)..(4,0+4);
\draw (-1,-1+4)..controls (1,-2+4) and (2,-2+4)..(4,0+4);

\draw (-2,0+4)--(0,0+4);
\draw[blue,thick] (-1,-1+4)--(0,0+4);

\draw (-2,0)--(-2,4);
\draw[blue,thick] (-1,-1)--(-1,3);
\draw[dashed,blue,thick] (0,0)--(0,4);
\draw[dashed,blue,thick] (3,0)--(3,4);
\draw[blue,thick] (4,0)--(4,4);

\draw[dashed,red,thick] (3,0)--(4,0);
\draw[dashed] (1.5,0.75)--(2,1.43);
\draw[dashed] (1.5,0.75)--(1.5,0.75+4);
\draw[dashed] (2,1.43)--(2,1.43+4);
\draw (1.5,0.75+4)--(2,1.43+4);
\draw[blue,thick] (3,4)--(4,4);

\filldraw (0,0) circle(2pt);
\filldraw (3,0) circle(2pt);
\filldraw (1.5,0.75) circle(2pt);

\node[above left] at (0,0) {$A_1$};
\node[above right] at (3,0) {$A_2$};
\node[above left] at (1.5,0.75) {$A_n$};

\node at (2.5,0.75) {\rotatebox{150}{$......$}};
\node at (2.5,0.75+4) {\rotatebox{150}{$......$}};

\begin{scope}[shift={(8,0)}]
\draw[dashed,red,thick] (0,0)..controls (1,1) and (2,1)..(3,0);
\draw[dashed] (0,0)..controls (1,-1) and (2,-1)..(3,0);

\draw[dashed,red,thick] (-2,0)..controls (1,2) and (2,2)..(4,0);
\draw (-1,-1)..controls (1,-2) and (2,-2)..(4,0);

\draw[red,thick] (-2,0)--(-1.1,0);
\draw[dashed,red,thick] (-0.9,0)--(0,0);
\draw[dashed] (-1,-1)--(0,0);

\draw[blue,thick] (0,0+4)..controls (1,1+4) and (2,1+4)..(3,0+4);
\draw (0,0+4)..controls (1,-1+4) and (2,-1+4)..(3,0+4);

\draw[blue,thick] (-2,0+4)..controls (1,2+4) and (2,2+4)..(4,0+4);
\draw (-1,-1+4)..controls (1,-2+4) and (2,-2+4)..(4,0+4);

\draw[blue,thick] (-2,0+4)--(0,0+4);
\draw (-1,-1+4)--(0,0+4);

\draw[blue,thick] (-2,0)--(-2,4);
\draw (-1,-1)--(-1,3);
\draw[dashed,blue,thick] (0,0)--(0,4);
\draw[dashed,blue,thick] (3,0)--(3,4);
\draw[blue,thick] (4,0)--(4,4);

\draw[dashed,red,thick] (3,0)--(4,0);
\draw[dashed,red,thick] (1.5,0.75)--(2,1.43);
\draw[dashed,blue,thick] (1.5,0.75)--(1.5,0.75+4);
\draw[dashed,blue,thick] (2,1.43)--(2,1.43+4);
\draw[blue,thick](1.5,0.75+4)--(2,1.43+4);
\draw[blue,thick] (3,4)--(4,4);

\filldraw (0,0) circle(2pt);
\filldraw (3,0) circle(2pt);
\filldraw (1.5,0.75) circle(2pt);

\node[above left] at (0,0) {$A_1$};
\node[above right] at (3,0) {$A_2$};
\node[above left] at (1.5,0.75) {$A_n$};

\node at (2.5,0.75) {\rotatebox{150}{$......$}};
\node at (2.5,0.75+4) {\rotatebox{150}{$......$}};
    
\end{scope}
\end{tikzpicture}
    \caption{$J\times I$ for Nontrivial $T$-Bundles}
    \label{fig:jtimesinontrivial}
\end{figure}

\begin{figure}
    \centering
\begin{tikzpicture}[scale=0.8]
\draw[dashed,green,thick] (0,0)..controls (1,1) and (2,1)..(3,0);
\draw[dashed] (0,0)..controls (1,-1) and (2,-1)..(3,0);

\draw[dashed] (-2,0)..controls (1,2) and (2,2)..(4,0);
\draw (-1,-1)..controls (1,-2) and (2,-2)..(4,0);

\draw[green,thick] (-2,0)--(-1.1,0);
\draw[green,thick] (-0.9,0)--(0,0);
\draw[dashed,green,thick] (-1,-1)--(0,0);

\draw[green,thick] (0,0+4)..controls (1,1+4) and (2,1+4)..(3,0+4);
%\draw (0,0+4)..controls (1,-1+4) and (2,-1+4)..(3,0+4);

\draw (-2,0+4)..controls (1,2+4) and (2,2+4)..(4,0+4);
\draw (-1,-1+4)..controls (1,-2+4) and (2,-2+4)..(4,0+4);

\draw (-2,0+4)--(0,0+4);
\draw (-1,-1+4)--(0,0+4);

\draw (-2,0)--(-2,4);
\draw (-1,-1)--(-1,3);
\draw[dashed,green,thick] (0,0)--(0,4);
%\draw[dashed,blue,thick] (3,0)--(3,4);
\draw[green,thick] (4,0)--(4,4);

\draw[dashed] (3,0)--(4,0);
\draw[dashed] (1.5,0.75)--(2,1.43);
\draw[dashed,green,thick] (1.5,0.75)--(1.5,0.75+4);
\draw[dashed,green,thick] (2,1.43)--(2,1.43+4);
\draw[green,thick](1.5,0.75+4)--(2,1.43+4);
\draw[green,thick] (3,4)--(4,4);

\filldraw (0,0) circle(2pt);
%\filldraw (3,0) circle(2pt);
\filldraw (1.5,0.75) circle(2pt);
\filldraw (1.5,0.75+4) circle(2pt);

\node[above left] at (0,0) {$A_1$};
\node[above right] at (3,0) {$A_2$};
\node[above left] at (1.5,0.75) {$A_n$};
\node[above left] at (1.5,0.75+4) {$A_n^{\prime}$};

\draw[dashed,green,thick](1.5,0.75)--(1.5,0.75-0.5);
\draw[dashed,green,thick] (2,1.43)--(2,1.43-0.5);

\node at (2.5,0.75) {\rotatebox{150}{$......$}};
\node at (2.5,0.75+4) {\rotatebox{150}{$......$}};

\draw [green,thick] (4,0)..controls (4,-1)..(5,-2);
\draw [dashed,green,thick] (3,0)--(3,-0.7);
\draw [green,thick] (3,-1)--(3,-2);
\draw [green,thick] (-2,0)--(-3,0);
\draw [green,thick] (-1,-1)--(-2,-2);

\node at (1.5,-3) {$n\geq2$};

\begin{scope}[shift={(8,0)}]
\draw[dashed] (0,0)..controls (1,1) and (2,1)..(3,0);
\draw[dashed] (0,0)..controls (1,-1) and (2,-1)..(3,0);

\draw[dashed] (-2,0)..controls (1,2) and (2,2)..(4,0);
\draw (-1,-1)..controls (1,-2) and (2,-2)..(4,0);

\draw[green,thick] (-2,0)--(-1.1,0);
\draw[green,thick] (-0.9,0)--(0,0);
\draw[dashed,green,thick] (-1,-1)--(0,0);

\draw (-2,0+4)..controls (1,2+4) and (2,2+4)..(4,0+4);
\draw (-1,-1+4)..controls (1,-2+4) and (2,-2+4)..(4,0+4);

\draw (-2,0+4)--(0,0+4);
\draw (-1,-1+4)--(0,0+4);

\draw (-2,0)--(-2,4);
\draw (-1,-1)--(-1,3);
\draw[dashed] (0,0)--(0,4);
\draw (4,0)--(4,4);

\node[above left] at (0,0) {$A_1$};

\draw [green,thick] (-2,0)--(-3,0);
\draw [green,thick] (-1,-1)--(-2,-2);

\draw[dashed] (3,0)--(3,-2);

\node at (1.5,-3) {$n=1$};
    
\end{scope}

\end{tikzpicture}
    \caption{$F_2$ for Nontrivial $T$-Bundles}    \label{fig:F2nontrivial}
\end{figure}

As shown in Figure~\ref{fig:jtimesinontrivial}, $J\times I$ consists of $n$ 3-balls whose 1-faces are indicated by red and blue lines. Similarly to the case of trivial bundles, we collapse $E$ to get $F_2$ as follows. First, we collapse the 3-ball indicated by the left diagram in Figure~\ref{fig:jtimesinontrivial} along the free 2-face enclosed by red lines. Next, we collapse the other 3-balls indicated by the right diagram in Figure~\ref{fig:jtimesinontrivial} by pushing upwards the 2-faces enclosed by red lines. We call this sequence of 3-deformations the \emph{nontrivial embedded disk move}.

When $n\geq2$, we get a new closed fake surface indicated by the left diagram in Figure~\ref{fig:F2nontrivial}. The original vertex $A_2$ is removed and new vertices $A_3^{\prime},\ldots,A_n^{\prime}$ are introduced. The original vertex $A_1$ remains, unlike in the case of trivial bundles, so the new surface $F_2$ contains $|\mathcal{S}_2(F_1)|+n-3$ vertices. 

In particular, for $n=2$, the number of vertices decreases by 1. Green edges are still connected together, so $F_2$ has the same number of connected components of 1-skeleton as $F_1$. In this case, the move is simply the $U^{-1}$ move from \cite{matveev2007algorithmic}. Note that for $n=3$, the number of vertices is unchanged. 

When $n=1$, we get a closed surface $F_2$ indicated by the right diagram in Figure~\ref{fig:F2nontrivial}. The original green circle $\partial D\subset\mathcal{S}_2(F_1)$ is changed to a dashed circle which represents a circle in a region. The other two remaining green edges are connected. Thus $F_2$ has the same number of connected components of 1-skeleton as $F_1$. The unique original vertex $A_1$ is removed, and so $F_2$ contains one fewer vertex than $F_1$.

\begin{figure}
    \centering
\begin{tikzpicture}[scale=1.5, thick]
      % First image (left)
  \begin{scope}
  % Define coordinates
  \coordinate (A) at (0,0);
  \coordinate (B) at (4,0);

  % Fill the region between the two curved edges in grey
  \fill[blue!30] (A) to[bend left=40] (B) to[bend left=40] (A);

  % Define coordinates along the top arc for the 4 vertices
  \path[bend left=40] (A) to coordinate[pos=0.2] (V1)
                           coordinate[pos=0.4] (V2)
                           coordinate[pos=0.6] (V3)
                           coordinate[pos=0.8] (V4) (B);

  % Draw curved arcs
  \draw[bend left=40] (A) to (B);
  \draw[red, very thick, bend right=40] (A) to (B);

  % Draw edges from each vertex (upward and downward)
  \draw (A) -- ++(0,0.8);
  \draw (A) -- ++(0,-0.8);
  \draw[red, very thick] (B) -- ++(0,0.8);
  \draw (B) -- ++(0,-0.8);
  \draw (V1) -- ++(0,0.8);
  \draw (V1) -- ++(0,-0.8);
  \draw (V2) -- ++(0,0.8);
  \draw (V2) -- ++(0,-0.8);
  \draw (V3) -- ++(0,0.8);
  \draw (V3) -- ++(0,-0.8);
  \draw (V4) -- ++(0,0.8);
  \draw (V4) -- ++(0,-0.8);

  % Draw the vertex circles (after edges so they appear on top)
  \fill (V1) circle (2pt);
  \fill (V2) circle (2pt);
  \fill (V3) circle (2pt);
  \fill (V4) circle (2pt);

    % Draw vertices
  \fill (A) circle (2pt) node[below] {};
  \fill (B) circle (2pt) node[below] {};

  % Optionally, add labels to the arcs
  \node at (2,0.9) {};
  \node at (2,-0.9) {};

  % a small blue square at (0,-2) to be used as a legend
  % add legend here
    %\fill[green!30] (0,-2) rectangle (0.5,-1.5);
    %\node at (1.9,-1.75) {Trivial $T$-bundle};

    % \fill[blue!30] (4,-2) rectangle (4.5,-1.5);
    % \node at (6.25,-1.75) {Nontrivial $T$-bundle};

  \end{scope}

    \draw[->, thick] (4.5,0) -- (5.5,0);

  % Second image (right) - copy of the first
  \begin{scope}[xshift=6cm]
  % Define coordinates
  \coordinate (A) at (0,0);
  \coordinate (B) at (4,0);

   % Define bottom points below A and B
  \coordinate (BottomA) at (0,-0.8);
  \coordinate (BottomB) at (4,-0.8);

  % Fill the region: top edge (A to B), down from B, across bottom, up to A
%   \fill[gray!30] (A) to[bend left=40] (B) -- (BottomB) -- (BottomA) -- cycle;

  % Define coordinates along the top arc for the 4 vertices
  \path[bend left=40] (A) to coordinate[pos=0.2] (V1)
                           coordinate[pos=0.4] (V2)
                           coordinate[pos=0.6] (V3)
                           coordinate[pos=0.8] (V4) (B);

  % Define top vertices
  \coordinate (T1) at (V1 |- 0,1.5);
  \coordinate (T2) at (V2 |- 0,1.5);
  \coordinate (T3) at (V3 |- 0,1.5);
  \coordinate (T4) at (V4 |- 0,1.5);

    %   \fill[gray!30] (0,0.0) to[bend left=28.5] (3.25,0.5) -- (T4) -- (T3) -- (T2) -- (T1) 
    % -- ++(-0.2,0) [rounded corners] -- (0,0.0) -- cycle;

        \fill[green!30] (0.75,0.5) to[bend left=21.5] (3.25,0.5) -- (T4) -- (T3) -- (T2) -- (T1) -- cycle;

        \fill[blue!30] (0,0.0) to[bend left=4] (0.75,0.5) -- (T1) -- ++(-0.15,0) -- ++(-0.05,-0.05) -- (0,0.0) -- cycle;

  % Draw curved arcs
  \draw[bend left=40] (A) to (B);

      % Draw edges from each vertex (upward and downward)
  \draw (A) -- ++(0,-0.8);
  \draw (A) -- ++(0,0.8);
  \draw (B) -- ++(0,-0.8);
  \draw (V1) -- (V1 |- 0,1.5);
  \draw (V1) -- ++(0,-0.8);
  \draw (V2) -- (V2 |- 0,1.5);
  \draw (V2) -- ++(0,-0.8);
  \draw (V3) -- (V3 |- 0,1.5);
  \draw (V3) -- ++(0,-0.8);
  \draw (V4) -- (V4 |- 0,1.5);
  \draw (V4) -- ++(0,-0.8);

  % Fill region from A to V4, bounded by curve below and red edges above
  
  % Connect top vertices with horizontal lines
\draw[red, very thick, rounded corners] (T1) --++(-0.2,0) -- (A);
  \draw[red, very thick] (T4) -- (4,1.5);
  \draw[red, very thick] (T1) -- (T2);
  \draw[red, very thick] (T2) -- (T3);
  \draw[red, very thick] (T3) -- (T4);

  % Draw the vertex circles (after edges so they appear on top)
  \fill (V1) circle (2pt);
  \fill (V2) circle (2pt);
  \fill (V3) circle (2pt);
  \fill (V4) circle (2pt);

    \draw (T1) -- ++(0,0.5);
    \draw (T2) -- ++(0,0.5);
    \draw (T3) -- ++(0,0.5);
    \draw (T4) -- ++(0,0.5);

      % Draw vertices
  \fill (A) circle (2pt) node[below] {};
%   \fill (B) circle (2pt) node[below] {};
  % Draw top vertex circles in red
  \fill[red] (T1) circle (2pt);
  \fill[red] (T2) circle (2pt);
  \fill[red] (T3) circle (2pt);
  \fill[red] (T4) circle (2pt);

  % Optionally, add labels to the arcs
  \node at (2,0.9) {};
  \node at (2,-0.9) {};
  \end{scope}
\end{tikzpicture}
\caption{The nontrivial embedded disk move. Embedded disks are colored, with green indicating a trivial $T$-bundle and blue indicating a nontrivial $T$-bundle.}
    \label{fig:nontrivial-embedded-disk-move}
\end{figure}
\begin{figure}
    \centering
   \begin{tikzpicture}[scale=1.5, thick]
      % First image (left)
  \begin{scope}
  % Define coordinates
  \coordinate (A) at (0,0);
  \coordinate (B) at (4,0);

  % Fill the region between the two curved edges in grey
  \fill[blue!30] (A) to[bend left=40] (B) to[bend left=40] (A);

  % Define coordinates along the top arc for the 4 vertices

  % Draw curved arcs
  \draw[bend left=40] (A) to (B);
  \draw[red, very thick, bend right=40] (A) to (B);

  % Draw edges from each vertex (upward and downward)
  \draw (A) -- ++(0,0.8);
  \draw (A) -- ++(0,-0.8);
  \draw[red, very thick] (B) -- ++(0,0.8);
  \draw (B) -- ++(0,-0.8);

    % Draw vertices
  \fill (A) circle (2pt) node[below] {};
  \fill (B) circle (2pt) node[below] {};

  % a small blue square at (0,-2) to be used as a legend
  % legend goes here
    % \fill[green!30] (0,-2) rectangle (0.5,-1.5);
    % \node at (1.9,-1.75) {Trivial $T$-bundle};

    % \fill[blue!30] (4,-2) rectangle (4.5,-1.5);
    % \node at (6.25,-1.75) {Nontrivial $T$-bundle};

  \end{scope}

    \draw[->, thick] (4.5,0) -- (5.5,0);
    % add node above arrow with text "T^{-1}
    \node at (5,0.3) {$U^{-1}$};

  % Second image (right) - copy of the first
  \begin{scope}[xshift=6cm]
  % Define coordinates
  \coordinate (A) at (0,0);
  \coordinate (B) at (4,0);

   % Define bottom points below A and B
  \coordinate (BottomA) at (0,-0.8);
  \coordinate (BottomB) at (4,-0.8);

  % Fill the region: top edge (A to B), down from B, across bottom, up to A
%   \fill[gray!30] (A) to[bend left=40] (B) -- (BottomB) -- (BottomA) -- cycle;

  % Define coordinates along the top arc for the 4 vertices
  \path[bend left=40] (A) to (B);

  % Define top vertices

    %   \fill[gray!30] (0,0.0) to[bend left=28.5] (3.25,0.5) -- (T4) -- (T3) -- (T2) -- (T1) 
    % -- ++(-0.2,0) [rounded corners] -- (0,0.0) -- cycle;

        % \fill[green!30] (0.75,0.5) to[bend left=21.5] (3.25,0.5) -- (T4) -- (T3) -- (T2) -- (T1) -- cycle;

        % \fill[blue!30] (0,0.0) to[bend left=20] (V1) -- (T1) -- ++(-0.95,0) -- ++(-0.05,-0.05) -- (0,0.0) -- cycle;

  % Draw curved arcs
  \draw[bend left=40] (A) to (B);

      % Draw edges from each vertex (upward and downward)
  \draw (A) -- ++(0,-0.8);
  \draw (A) -- ++(0,0.8);
  \draw (B) -- ++(0,-0.8);

  % Connect top vertices with horizontal lines
\draw[red, very thick, rounded corners] (T1) --++(-0.2,0) -- (A);
  \draw[red, very thick] (T1) -- (4,1.5);

  % Draw the vertex circles (after edges so they appear on top)
%   \fill (V1) circle (2pt);

    % \draw (T1) -- ++(0,0.5);

      % Draw vertices
  \fill (A) circle (2pt) node[below] {};
%   \fill (B) circle (2pt) node[below] {};
  % Draw top vertex circles in red
%   \fill[red] (T1) circle (2pt);
  \end{scope}
\end{tikzpicture}
    \caption{The nontrivial embedded disk move in the case $n=2$. This is the maximal case in which the complexity is reduced. In this case, the move is equivalent to the $U^{-1}$ move from \cite{matveev2007algorithmic}.}
    \label{fig:placeholder}
\end{figure}

\bigskip

In the remainder of the paper, we refer to these two moves the \emph{trivial embedded disk move} and the \emph{nontrivial embedded disk move}. When we apply these moves to an embedded disk, we say that we do so ``along" an edge, which refer to the edge(s) of the embedded disk that are removed. 

\section{Reducing the Complexity of a Fake Surface}
\label{sec:crl}
In this section, we prove the following theorem, which is the main result of this work:
\begin{theorem}[Precise statement of Theorem~\ref{thm:main}]
\label{theorem:ac}

For any surface $F\in\mathcal{F}_C(1,t)$ where $t\leq 5$, there exists a zigzag
$$F\swarrow F_1\searrow F_2\swarrow...\searrow F_{n-1}\swarrow F_n\searrow *
$$
of elementary collapses, each either leftwards or rightwards, such that $\text{dim}(F_i)\leq3$.
\end{theorem}

\subsection{Complexity Reduction Lemma}

We begin by proving a lemma consisting of conditions under which a fake surface can be reduced in complexity.

\begin{lemma}[Complexity Reduction Lemma]\label{lem:crl} A fake surface $F$ can be 3-deformed to another fake surface of lower complexity if $F$ satisfies one of the following criteria:
\begin{itemize}
\item[1)] $F$ contains an embedded disk with at most three vertices and a trivial $T$-bundle.
\item[2)] $F$ contains an embedded disk with at most two vertices and a nontrivial $T$-bundle.
\item[3)] $F$ contains two embedded disks with three vertices that share exactly one vertex.
\item[4)] $F$ contains one embedded disk with 3 vertices and another embedded disk with 3 or 4 vertices, and the two embedded disks share exactly one edge. 
\item[5)] $F$ contains an embedded disk with 4 vertices and trivial $T$-bundle and another disk with 4 vertices, and either the two disks share exactly two non-adjacent edges or the embedded disk intersects the other disk in one edge twice. If the second disk has a trivial $T$-bundle, then the condition that two edges are non-adjacent can be removed. 
\item[6)] $F$ contains an embedded disk with 4 vertices and trivial $T$-bundle, another disk with 4 vertices that intersects the embedded disk in two consecutive edges, and another disk with 4 or 5 vertices as in Figures~\ref{fig:4+4+51} and \ref{fig:4+4+52}.
\item[7)] $F$ contains an embedded disk with 5 vertices and trivial $T$-bundle and another embedded disk with 4 vertices and trivial $T$-bundle, and the two disks share exactly two non-adjacent edges. 
\item[8)] $F$ contains an embedded disk with 5 vertices and trivial $T$-bundle and three disks with 4 vertices as in either Figures~\ref{fig:5+4+4+41} or \ref{fig:5+4+4+42}.

\end{itemize}
    
\end{lemma}

Our proof of the complexity reduction lemma goes case by case as follows.

\subsection*{Cases 1 and 2}
Cases 1 and 2 follow directly from applying the trivial embedded disk move and nontrivial embedded disk move described in Section~\ref{sec:embedded-disk-moves}.
\subsection*{Case 3}

It suffices to consider the case where both embedded disks have nontrivial $T$-bundles, since otherwise we can reduce the complexity by Case 1. In Figure~\ref{fig:3+3 vertex}, let the $T$-bundle of the red disk $[7,-12,4]$ be given by 
\[
(3,7,15),(-15,-12,13),(-13,4,1),(-1,7,16),(-16,-12,14),(-14,4,-3),
\] 
and the $T$-bundle of the blue disk $[1,2,3]$ be given by
\[
(4,1,5),(-5,2,6),(-6,3,7),(-7,1,8),(-8,2,9),(-9,3,-4).
\]

We now apply the nontrivial embedded disk move to the blue disk along edge $2$. After applying the move, the red disk becomes another embedded disk with 3 vertices whose $T$-bundle given by 
\[
(11,7,15),(-15,-12,13),(-13,4,-11), 
\]
where edge 11 is the new diagonal edge as in Figure~\ref{fig:nontrivial-embedded-disk-move}. This is a trivial $T$-bundle, so we can reduce the complexity by Case 1.

\begin{figure}[ht]
    \centering
    \caption{Configuration for Case 3}
\begin{tikzpicture}[scale=0.8]
\draw (-1,0)..controls (1,-2) and (3,-2)..(5,0);
\draw (-1,0)..controls (1,2) and (3,2)..(5,0);

\draw (3,0.6-1)--(3,0.6+2);
\draw[rounded corners] (-1,0)--(-1,2.5)--(-2,2.5);
\draw[rounded corners] (-1,0)--(-1,-2.5)--(-2,-2.5);
\draw (-2,2.5)--(-2,-2.5);
\draw (-2,2.5)--(-2,3.5);
\draw (-2,2.5)--(-3,2.5);
\draw (-2,-2.5)--(-2,-3.5);
\draw (-2,-2.5)--(-3,-2.5);

\draw (5,-2.5)--(5,2.5);

\filldraw (3,1.35) circle (2pt);
\filldraw (-1,0) circle (2pt);
\filldraw (5,0) circle (2pt);
\filldraw (-2,2.5) circle (2pt);
\filldraw (-2,-2.5) circle (2pt);

\node at (2,-2){$1$};
\node at (4,1.5){$2$};
\node at (-0.5,2){$4$};
\node at (5,3){$5$};
\node at (3,3){$6$};
\node at (-0.5,-2){$7$};
\node at (5,-3){$8$};
\node at (3,-1){$9$};
\node at (-2.5,0){$12$};
\node at (-2,4){$13$};
\node at (-3.5,2.5){$14$};
\node at (-2,-4){$15$};
\node at (-3.5,-2.5){$16$};

% Replace all arrow nodes with proper Stealth arrows
% > at (3,2) rotated 90° (pointing up)
\draw[-{Stealth[length=6pt]}] (3,2) -- ++(90:0.01);
% > at (-1,2) rotated -90° (pointing down)
\draw[-{Stealth[length=6pt]}] (-1,2) -- ++(-90:0.01);
% > at (5,2) rotated 90° (pointing up)
\draw[-{Stealth[length=6pt]}] (5,2) -- ++(90:0.01);
% > at (3,0.5) rotated -90° (pointing down)
\draw[-{Stealth[length=6pt]}] (3,0.5) -- ++(-90:0.01);
% > at (-1,-2) rotated -90° (pointing down)
\draw[-{Stealth[length=6pt]}] (-1,-2) -- ++(-90:0.01);
% > at (5,-2) rotated -90° (pointing down)
\draw[-{Stealth[length=6pt]}] (5,-2) -- ++(-90:0.01);
% > at (2,-1.5) not rotated (pointing right)
\draw[-{Stealth[length=6pt]}] (2,-1.5) -- ++(0:0.01);
\node at (2,2){$3$};
% < at (-2,0) rotated 90° (pointing down)
\draw[-{Stealth[length=6pt]}] (-2,0) -- ++(-90:0.01);
% < at (2,1.5) not rotated (pointing left)
\draw[-{Stealth[length=6pt]}] (2,1.5) -- ++(180:0.01);
% < at (3.7,1) rotated -45° (pointing down-left at 135°)
\draw[-{Stealth[length=6pt]}] (3.7,1.04) -- ++(140:0.01);

% > at (-2,3) rotated 90° (pointing up)
\draw[-{Stealth[length=6pt]}] (-2,3) -- ++(90:0.01);
% < at (-2.5,2.5) not rotated (pointing left)
\draw[-{Stealth[length=6pt]}] (-2.5,2.5) -- ++(180:0.01);
% > at (-2,-3) rotated 270° (pointing down)
\draw[-{Stealth[length=6pt]}] (-2,-3) -- ++(270:0.01);
% < at (-2.5,-2.5) not rotated (pointing left)
\draw[-{Stealth[length=6pt]}] (-2.5,-2.5) -- ++(180:0.01);

\filldraw[opacity=0.3,blue] (-1,0)..controls (1,-2) and (3,-2)..(5,0)..controls (3,2) and (1,2)..(-1,0);
\filldraw[opacity=0.3,red,rounded corners] (-1,0)--(-1,2.5)--(-2,2.5)--(-2,-2.5)--(-1,-2.5)--(-1,0);

\end{tikzpicture}
    \caption*{Two embedded disks with 3 vertices that share exactly one vertex.}
    \label{fig:3+3 vertex}
\end{figure}

\subsection*{Case 4}
It suffices to assume that the $T$-bundle of the blue disk in Figure~\ref{fig:3+3 3+4 edge} is nontrivial, since otherwise we can apply case 1. We apply the nontrivial embedded disk move to the blue disk along the shared edge such that the disks are changed to two disks that share exactly one vertex. The two possibilities are shown in Figure~\ref{fig:3+3 3+4 edge}. In each case, we can apply Case 1, 2, or 3 to reduce the complexity. 

\begin{figure}[ht]
    \centering
    \caption{Configurations for Case 4}
\begin{tikzpicture}
\draw (0,0)--(1,0)--(1,1)--(0,0);
\draw (1,0)--(2,1)--(1,1);

\draw[->,line width=1pt] (3,0.5)--(4,0.5);

\draw (5,0)--(6,0)--(6,1)--(5,0);
\draw (6,1)..controls(6.5,1.5)..(7,1)..controls(6.5,0.5)..(6,1);

\filldraw[opacity=0.3,blue] (0,0)--(1,0)--(1,1)--(0,0);
\filldraw[opacity=0.3,red] (1,0)--(2,1)--(1,1)--(1,0);
\filldraw[opacity=0.3,blue] (5,0)--(6,0)--(6,1)--(5,0);
\filldraw[opacity=0.3,red] (6,1)..controls(6.5,1.5)..(7,1)..controls(6.5,0.5)..(6,1);

\filldraw (0,0) circle(2pt);
\filldraw (1,0) circle(2pt);
\filldraw (1,1) circle(2pt);
\filldraw (2,1) circle(2pt);
\filldraw (5,0) circle(2pt);
\filldraw (6,0) circle(2pt);
\filldraw (6,1) circle(2pt);
\filldraw (7,1,0) circle(2pt);

\draw (0,-3)--(1,-3)--(1,-2)--(0,-3);
\draw (1,-3)--(2,-3)--(2,-2)--(1,-2);

\draw[->,line width=1pt] (3,-2.5)--(4,-2.5);

\draw (5,-2)--(5,-3)--(6,-3)--(5,-2);
\draw (6,-3)--(7,-3)--(7,-2)--(6,-3);

\filldraw[opacity=0.3,blue] (0,-3)--(1,-3)--(1,-2)--(0,-3);
\filldraw[opacity=0.3,red] (1,-3)--(2,-3)--(2,-2)--(1,-2)--(1,-3);
\filldraw[opacity=0.3,blue] (5,-2)--(5,-3)--(6,-3)--(5,-2);
\filldraw[opacity=0.3,red] (6,-3)--(7,-3)--(7,-2)--(6,-3);

\filldraw (0,-3) circle(2pt);
\filldraw (1,-3) circle(2pt);
\filldraw (1,-2) circle(2pt);
\filldraw (2,-3) circle(2pt);
\filldraw (2,-2) circle(2pt);
\filldraw (5,-3) circle(2pt);
\filldraw (5,-2) circle(2pt);
\filldraw (6,-3) circle(2pt);
\filldraw (7,-3) circle(2pt);
\filldraw (7,-2) circle(2pt);

\end{tikzpicture}
\vspace{0.8cm}
    \caption*{\textbf{Left:} An embedded disk with 3 vertices and an embedded disk with 3 or 4 vertices that share exactly one edge. \textbf{Right:} The resulting configuration after applying the nontrivial embedded disk move to the blue disk.}
    \label{fig:3+3 3+4 edge}
\end{figure}
\ 

\subsection*{Case 5}

The two possible configurations for Case 5 are shown in Figure~\ref{fig:4+4}. We apply the trivial embedded disk move to the blue disk along the thickened edges. The red disk will be changed to a disk with 2 vertices, allowing us to reduce the complexity using Case 1 or 2. If the red disk has a trivial $T$-bundle, the two edges on which the disks intersect can be adjacent. This is because the trivial embedded disk move will not change the triviality of disks, and in the case of adjacent edges, the second disk will become a disk with no more than 3 vertices and trivial $T$-bundle. Thus we can reduce the complexity using Case 1. 

\begin{figure}[ht]
    \centering
        \caption{Configurations for Case 5}
\begin{tikzpicture}
\draw (0,0)--(3,0);
\draw[rounded corners](3,0)--(3,-1)--(0,-1)--(0,0);
\draw[rounded corners] (0,0)--(0,2)--(3,2)--(3,0);
\draw[rounded corners] (1,0)--(1,1)--(2,1)--(2,0);
\draw[line width=2pt] (0,0)--(1,0);
\draw[line width=2pt] (2,0)--(3,0);

\filldraw (0,0) circle(2pt);
\filldraw (1,0) circle(2pt);
\filldraw (2,0) circle(2pt);
\filldraw (3,0) circle(2pt);

\filldraw[opacity=0.3,blue,rounded corners] (0,0)--(3,0)--(3,-1)--(0,-1)--(0,0);
\filldraw[opacity=0.3,red,rounded corners] (0,0)--(1,0)--(1,1)--(2,1)--(2,0)--(3,0)--(3,2)--(0,2)--(0,0);

\draw (0+5,0)--(3+5,0);
\draw[rounded corners](3+5,0)--(3+5,-1)--(0+5,-1)--(0+5,0);

\draw (7,0)..controls(7.5,1)..(8,0)..controls(7.5,-1)..(7,0);

\filldraw (0+5,0) circle(2pt);
\filldraw (1+5,0) circle(2pt);
\filldraw (2+5,0) circle(2pt);
\filldraw (3+5,0) circle(2pt);

\draw[line width=2pt] (5,0)--(6,0);
\draw[line width=2pt] (7,0)--(8,0);

\node at (7.5,1){$1$};
% \draw[-{Stealth[length=6pt]}] (7.2,0.3) -- ++(0.01,0.02);
\node at (7.5,0.25){$2$};
\node at (7.5,-0.5){$3$};

\filldraw[opacity=0.3,blue,rounded corners] (5,0)--(8,0)--(8,-1)--(5,-1)--(5,0);
\filldraw[opacity=0.3,red,rounded corners](7,0)..controls(7.5,1)..(8,0)..controls(7.5,-1)..(7,0);

\end{tikzpicture}
\vspace{0.8cm}
\caption*{\textbf{Left:} Two embedded disks with 4 vertices that share two non-adjacent edges. \textbf{Right:} Two embedded disks with 4 vertices where the red disk intersects the blue disk at one edge twice. The red disk's attaching map is $[2,1,2,3]$.}
    \label{fig:4+4}
\end{figure}

\subsection*{Case 6}
The two possible configurations in Case 6 are given in Figures~\ref{fig:4+4+51} and \ref{fig:4+4+52}. It is sufficient to assume that the $T$-bundles of the second and third disks are nontrivial, since otherwise we can reduce the complexity by Case 5. In the first case, shown in Figure~\ref{fig:4+4+51}, we do the trivial embedded disk move with the blue disk $[1,2,3,4]$ along edges 2 and 4. We get an embedded disk with 3 vertices and an embedded disk with 3 or 4 vertices that share exactly one edge, and the resulting fake surface has the same complexity as before. Thus, we can apply Case 4. In the second case, shown in Figure~\ref{fig:4+4+52}, we again do the trivial embedded disk move with the blue disk $[1,2,3,4]$ along edges 2 and 4. This again gives an embedded disk with 3 vertices and an embedded disk with 3 or 4 vertices that share exactly one edge, and again the resulting fake surface has the same complexity as before. We apply Case 4 to reduce the complexity.  

\begin{figure}
    \centering
    \caption{First Configuration for Case 6}
\begin{tikzpicture}
\draw (0,0)--(2,0)--(2,2)--(0,2)--(0,0);
\draw (2,0)..controls(3,1)..(2,2);
\draw (0,2)..controls(1,3)..(2,2);
\draw (0,2)..controls(-1,1)..(0,0);

\filldraw (0,0) circle(2pt);
\filldraw (2,0) circle(2pt);
\filldraw (2,2) circle(2pt);
\filldraw (0,2) circle(2pt);

\filldraw[opacity=0.3,blue,rounded corners] (0,0)--(2,0)--(2,2)--(0,2)--(0,0);
\filldraw[opacity=0.3,red,rounded corners] (2,0)--(2,2)..controls(3,1)..(2,0);
\filldraw[opacity=0.3,red,rounded corners] (0,2)--(2,2)..controls(1,3)..(0,2);
\filldraw[opacity=0.3,green,rounded corners] (0,0)--(0,2)..controls(-1,1)..(0,0);
\filldraw[opacity=0.3,green,rounded corners] (0,2)--(2,2)..controls(1,3)..(0,2);

\node at (1,-0.3){$1$};
\node at (2.3,1){$2$};
\node at (1,2.3){$3$};
\node at (-0.3,1){$4$};
\node at (3,1){$6$};
\node at (1,3){$5$};
\node at (-1.7,1){$7$ or $7,8$};

% Replace all arrow nodes in left part with proper Stealth arrows
% > at (1,0) not rotated (pointing right)
\draw[-{Stealth[length=6pt]}] (1.1,0) -- ++(0:0.01);
% > at (2,1) rotated 90° (pointing up)
\draw[-{Stealth[length=6pt]}] (2,1.1) -- ++(90:0.01);
% < at (1,2) not rotated (pointing left)
\draw[-{Stealth[length=6pt]}] (0.9,2) -- ++(180:0.01);
% > at (0,1) rotated 270° (pointing down)
\draw[-{Stealth[length=6pt]}] (0,0.9) -- ++(270:0.01);
% < at (1,2.75) not rotated (pointing left)
\draw[-{Stealth[length=6pt]}] (0.97,2.75) -- ++(180:0.01);
% < at (2.75,1) rotated 90° (pointing down)
\draw[-{Stealth[length=6pt]}] (2.75,0.97) -- ++(-90:0.01);
% > at (-0.75,1) rotated 270° (pointing down)
\draw[-{Stealth[length=6pt]}] (-0.75,0.97) -- ++(270:0.01);

\draw[->,line width=1pt] (3.5,1)--(4.5,1);

\draw (6+1,2)--(7+1,2)--(7+1,0)--(6+1,0)--(6+1,2);
\draw (7+1,2)--(8+1,1)--(7+1,0);

\filldraw (7+1,2) circle(2pt);
\filldraw (7+1,0) circle(2pt);
\filldraw (6+1,0) circle(2pt);
\filldraw (8+1,1) circle(2pt);

\filldraw[opacity=0.3,green] (6+1,2)--(7+1,2)--(7+1,0)--(6+1,0)--(6+1,2);
\filldraw[opacity=0.3,red] (7+1,2)--(8+1,1)--(7+1,0)--(7+1,2);

\node at (6.7+1,1){$5$};
\node at (6.5+1,-0.2){$1'$};
\node at (7.5+1,2){$1$};
\node at (7.5+1,0){$6$};
\node at (5+1,1){$7$ or $7,8$};

\end{tikzpicture}
    \caption*{\textbf{Left:} There is an embedded disk $[1,2,3,4]$, a disk $[2,5,-3,6]$, and a disk $[-4,-5,3,7]$ or $[-4,-5,3,7,8]$. \textbf{Right:} The resulting configuration has an embedded disk with 3 vertices and an embedded disk with 3 or 4 vertices that share exactly one edge.}
    \label{fig:4+4+51}
\end{figure}

\begin{figure}
    \centering
    \caption{Second Configuration for Case 6}
\begin{tikzpicture}
\draw (0,0)--(2,0)--(2,2)--(0,2)--(0,0);
\draw[rounded corners] (2,2)..controls(2,3)..(3,3)..controls(3,2)..(2,2);
\draw (0,0)..controls(-1,0.5)..(-2,0);
\draw (0,0)..controls(-1,-0.5)..(-2,0);
\draw (0,2)..controls(1.5,1.5)..(2,0);
\draw (0,2)..controls(0.5,0.5)..(2,0);

\filldraw (0,0) circle(2pt);
\filldraw (2,0) circle(2pt);
\filldraw (2,2) circle(2pt);
\filldraw (0,2) circle(2pt);

\filldraw[opacity=0.3,blue,rounded corners] (0,0)--(2,0)--(2,2)--(0,2)--(0,0);
\filldraw[opacity=0.3,red,rounded corners] (0,2)--(2,2)..controls(2,3)..(3,3)..controls(3,2)..(2,2)--(2,0)..controls(1.5,1.5)..(0,2);
\filldraw[opacity=0.3,green,rounded corners] (0,2)--(0,0)..controls(-1,0.5)..(-2,0)..controls(-1,-0.5)..(0,0)--(2,0)..controls(0.5,0.5)..(0,2);

\node at (1,-0.3){$1$};
\node at (2.3,1){$2$};
\node at (1,2.3){$3$};
\node at (-0.3,1){$4$};
\node at (3.3,2.5){$5$};
\node at (1.5,1.5){$6$};
\node at (0.4,0.4){$7$};
\node at (-1,-0.7){$8$ or $8,9$};

\draw[->,line width=1pt] (3.5,1)--(4.5,1);

\draw (6+1,2)--(7+1,2)--(7+1,0)--(6+1,0)--(6+1,2);
\draw (7+1,2)--(8+1,1)--(7+1,0);

\filldraw (7+1,2) circle(2pt);
\filldraw (7+1,0) circle(2pt);
\filldraw (6+1,0) circle(2pt);
\filldraw (8+1,1) circle(2pt);

\filldraw[opacity=0.3,green] (6+1,2)--(7+1,2)--(7+1,0)--(6+1,0)--(6+1,2);
\filldraw[opacity=0.3,red] (7+1,2)--(8+1,1)--(7+1,0)--(7+1,2);

\node at (6.7+1,1){$1'$};
\node at (6.5+1,-0.2){$7$};
\node at (7.5+1,2){$5$};
\node at (7.5+1,0){$6$};
\node at (5+1,1){$8$ or $8,9$};

% Replace all arrow nodes in left part with proper Stealth arrows
% > at (1,0) not rotated (pointing right)
\draw[-{Stealth[length=6pt]}] (1,0) -- ++(0:0.01);
% > at (2,1) rotated 90° (pointing up)
\draw[-{Stealth[length=6pt]}] (2,1) -- ++(90:0.01);
% < at (1,2) not rotated (pointing left)
\draw[-{Stealth[length=6pt]}] (1,2) -- ++(180:0.01);
% > at (0,1) rotated 270° (pointing down)
\draw[-{Stealth[length=6pt]}] (0,1) -- ++(270:0.01);
% > at (2.5,2) not rotated (pointing right)
\draw[-{Stealth[length=6pt]}] (2.5,2) -- ++(0:0.01);
% > at (1.35,1.35) rotated 315° (pointing down-right, which is -45°)
\draw[-{Stealth[length=6pt]}] (1.35,1.38) -- ++(-45:0.01);
% > at (0.65,0.65) rotated 135° (pointing up-left)
\draw[-{Stealth[length=6pt]}] (0.65,0.62) -- ++(135:0.01);
% < at (-1,0.4) not rotated (pointing left)
\draw[-{Stealth[length=6pt]}] (-1,0.38) -- ++(180:0.01);

\draw[->,line width=1pt] (3.5,1)--(4.5,1);

\draw (6+1,2)--(7+1,2)--(7+1,0)--(6+1,0)--(6+1,2);
\draw (7+1,2)--(8+1,1)--(7+1,0);

\end{tikzpicture}
    \caption*{\textbf{Left:} There is an embedded disk $[1,2,3,4]$, a disk $[2,5,3,6]$, and a disk $[1,7,4,8]$ or $[1,7,4,8,9]$. \textbf{Right:} The resulting configuration has an embedded disk with 3 vertices and an embedded disk with 3 or 4 vertices that share exactly one edge.}
    \label{fig:4+4+52}
\end{figure}
\ 

\subsection*{Case 7}

The configuration of Case 7 is shown in Figure~\ref{fig:5+4}. In this case, we do the trivial embedded disk move to the blue disk with 5 vertices along the thick edges as in Figure~\ref{fig:5+4}. This increases the complexity of the resulting fake surface by 1. However, the red disk is changed to a disk with two vertices and trivial $T$-bundle. We then apply the trivial embedded disk move to this disk, which decreases the complexity by 2. Thus, the resulting surface is of lower complexity than the original one. 

\begin{figure}
    \centering
    \caption{Configuration for Case 7}
\begin{tikzpicture}
\draw (0,0)--(3,0);
\draw[rounded corners](3,0)--(3,-1)--(0,-1)--(0,0);
\draw[rounded corners] (0,0)--(0,2)--(3,2)--(3,0);
\draw[rounded corners] (1,0)--(1,1)--(2,1)--(2,0);
\draw[line width=2pt] (0,0)--(1,0);
\draw[line width=2pt] (2,0)--(3,0);

\filldraw (0,0) circle(2pt);
\filldraw (1,0) circle(2pt);
\filldraw (2,0) circle(2pt);
\filldraw (3,0) circle(2pt);
\filldraw (1.5,-1) circle(2pt);

\filldraw[opacity=0.3,blue,rounded corners] (0,0)--(3,0)--(3,-1)--(0,-1)--(0,0);
\filldraw[opacity=0.3,red,rounded corners] (0,0)--(1,0)--(1,1)--(2,1)--(2,0)--(3,0)--(3,2)--(0,2)--(0,0);

\end{tikzpicture}
    \caption*{An embedded disk with 5 vertices and an embedded disk with 4 vertices that share two non-adjacent edges.}
    \label{fig:5+4}
\end{figure}

\subsection*{Case 8}

The first case is shown in Figure~\ref{fig:5+4+4+41}. We can assume that $[2,9,4,8]$ and $[10,-3,9,-3]$ have nontrivial $T$-bundles, since otherwise we can reduce the complexity by Case 7. We do the trivial embedded disk move to $[1,2,3,4,5]$ along edges 1 and 3, which increases the complexity by 1. The resulting fake surface has three embedded disks as in the right side of Figure~\ref{fig:5+4+4+41}, where the green disk with 2 vertices and the red disk with 5 vertices have nontrivial $T$-bundles. We then do the nontrivial embedded disk move to the green disk to remove edge 9 from the red disk and go back to the original complexity. The new red disk has 4 edges, and so we can apply Case 4 to the new red disk and the gray disk. We can use a similar method for the second case in Figure~\ref{fig:5+4+4+42}. 

\begin{figure}
   \centering
   \caption{First Configuration for Case 8}
\begin{tikzpicture}
\draw (0,0)--(-1,-1)--(-1,-2)--(1,-2)--(1,-1)--(0,0);
\draw (-1,-1)--(1,-1);
\draw[rounded corners] (0,0)..controls(-0.5,0.5)..(0,1)..controls(0.5,0.5)..(0,0);
\draw (-1,-1)..controls(0,-1.5)..(1,-1);
\draw (-1,-2)..controls(0,-2.5)..(1,-2);
\draw[rounded corners] (-1,-2)--(-1,-3)--(1,-3)--(1,-2);

\filldraw (0,0) circle(2pt);
\filldraw (-1,-1) circle(2pt);
\filldraw (1,-1) circle(2pt);
\filldraw (-1,-2) circle(2pt);
\filldraw (1,-2) circle(2pt);

\filldraw[opacity=0.3,blue,rounded corners] (0,0)--(-1,-1)--(-1,-2)--(1,-2)--(1,-1)--(0,0);
\filldraw[opacity=0.3,red,rounded corners] (-1,-1)..controls(0,-1.5)..(1,-1)--(1,-2)..controls(0,-2.5)..(-1,-2)--(-1,-1);
\filldraw[opacity=0.3,green,rounded corners] (-1,-2)--(1,-2)--(1,-3)--(-1,-3)--(-1,-2);
\filldraw[opacity=0.3,rounded corners] (-1,-1)--(0,0)..controls(-0.5,0.5)..(0,1)..controls(0.5,0.5)..(0,0)--(1,-1)--(-1,-1);

\node at (-1,-0.5){$1$};
\node at (-1.4,-1.5){$2$};
\node at (0.5,-1.8){$3$};
\node at (1.4,-1.5){$4$};
\node at (1,-0.5){$5$};
\node at (0.6,0.7){$6$};
\node at (0,-0.7){$7$};
\node at (-0.5,-1.45){$8$};
\node at (-0.5,-2.5){$9$};
\node at (0,-3.3){$10$};

\draw[->,line width=1pt] (2,-1.5)--(3,-1.5);

\draw (5.5,0)--(4.5,-1)--(4.5,-2)--(6.5,-2)--(6.5,-1)--(5.5,0);
\draw (4.5,-2)--(5.5,-3)--(6.5,-2);
\draw[rounded corners] (5.5,0)--(6.5,0)--(6.5,-1);

\filldraw (5.5,0) circle(2pt);
\filldraw (4.5,-1) circle(2pt);
\filldraw (4.5,-2) circle(2pt);
\filldraw (6.5,-2) circle(2pt);
\filldraw (6.5,-1) circle(2pt);
\filldraw (5.5,-3) circle(2pt);

\filldraw[opacity=0.3,red,rounded corners] (5.5,0)--(4.5,-1)--(4.5,-2)--(6.5,-2)--(6.5,-1)--(5.5,0);
\filldraw[opacity=0.3,green,rounded corners] (5.5,0)--(6.5,0)--(6.5,-1)--(5.5,0);
\filldraw[opacity=0.3,rounded corners] (4.5,-2)--(6.5,-2)--(5.5,-3)--(4.5,-2);

\node at (5.75,-0.75){$9$};
\node at (7,-0.5){$10$};
\node at (7,-1.5){$4$};
\node at (5.5,-1.7){$5$};
\node at (4,-1.5){$8$};
\node at (4.75,-0.25){$4'$};
\node at (6.25,-2.75){$7$};
\node at (4.75,-2.75){$6$};

% Replace all arrow nodes in left part with proper Stealth arrows
% < at (-0.5,-0.5) rotated 45° (pointing up-right at 225°)
\draw[-{Stealth[length=6pt]}] (-0.5,-0.5) -- ++(225:0.01);
% < at (-1,-1.5) rotated 90° (pointing down at -90°)
\draw[-{Stealth[length=6pt]}] (-1,-1.6) -- ++(-90:0.01);
% < at (0,-2) rotated 180° (pointing right at 0°)
\draw[-{Stealth[length=6pt]}] (0.1,-2) -- ++(0:0.01);
% < at (1,-1.5) rotated 270° (pointing up at 90°)
\draw[-{Stealth[length=6pt]}] (1,-1.4) -- ++(90:0.01);
% < at (0.5,-0.5) rotated 315° (pointing down-right at 135°)
\draw[-{Stealth[length=6pt]}] (0.5,-0.5) -- ++(135:0.01);
% < at (0,-2.35) rotated 180° (pointing right at 0°)
\draw[-{Stealth[length=6pt]}] (0.1,-2.37) -- ++(0:0.01);
% < at (0,-1.35) not rotated (pointing left at 180°)
\draw[-{Stealth[length=6pt]}] (-0.04,-1.36) -- ++(180:0.01);
% < at (0,-3) rotated 180° (pointing right at 0°)
\draw[-{Stealth[length=6pt]}] (0.1,-3) -- ++(0:0.01);
% > at (0,0.9) not rotated (pointing right at 0°)
\draw[-{Stealth[length=6pt]}] (0.05,0.95) -- ++(0:0.01);
% < at (0,-1) rotated 180° (pointing right at 0°)
\draw[-{Stealth[length=6pt]}] (0.1,-1) -- ++(0:0.01);

\end{tikzpicture}
    \caption*{\textbf{Left:} There are embedded disks [1,2,3,4,5], [1,7,5,6], [2,9,4,8], and [10,-3,9,-3]. \textbf{Right:} The resulting configuration has an embedded disk with 2 vertices, and embedded disk with 5 vertices, and an embedded disk with 3 vertices.}
   \label{fig:5+4+4+41}
\end{figure}

\begin{figure}
    \centering
    \caption{Second Configuration for Case 8}
\begin{tikzpicture}
\draw (0,0)--(-1,-1)--(-1,-2)--(1,-2)--(1,-1)--(0,0);
\draw (-1,-1)--(1,-1);
\draw (-1,-1)..controls(0,-1.5)..(1,-1);
\draw (-1,-2)..controls(0,-2.5)..(1,-2);
\draw[rounded corners] (0,0)--(2,0)--(2,-2)--(1,-2);
\draw[rounded corners] (0,0)--(-2,0)--(-2,-2)--(-1,-2);

\filldraw (0,0) circle(2pt);
\filldraw (-1,-1) circle(2pt);
\filldraw (1,-1) circle(2pt);
\filldraw (-1,-2) circle(2pt);
\filldraw (1,-2) circle(2pt);

\filldraw[opacity=0.3,blue,rounded corners] (0,0)--(-1,-1)--(-1,-2)--(1,-2)--(1,-1)--(0,0);
\filldraw[opacity=0.3,red,rounded corners] (-1,-1)..controls(0,-1.5)..(1,-1)--(1,-2)..controls(0,-2.5)..(-1,-2)--(-1,-1);
\filldraw[opacity=0.3,green,rounded corners] (0,0)--(-2,0)--(-2,-2)--(-1,-2)--(-1,-1)--(1,-1)--(0,0);
\filldraw[opacity=0.3,rounded corners] (0,0)--(2,0)--(2,-2)--(1,-2)--(1,-1)--(-1,-1)--(0,0);

\node at (-0.8,-0.5){$1$};
\node at (-1.3,-1.5){$2$};
\node at (0.3,-1.8){$3$};
\node at (1.3,-1.5){$4$};
\node at (0.8,-0.5){$5$};
\node at (2.3,-1){$6$};
\node at (0,-0.7){$7$};
\node at (-0.5,-1.45){$8$};
\node at (-0.5,-2.5){$9$};
\node at (-2.3,-1){$10$};

% Replace all arrow nodes with proper Stealth arrows
% Arrow at (-0.5,-0.5) pointing at 45 degrees (up-right)
\draw[-{Stealth[length=6pt]}] (-0.5,-0.5) -- ++(45:-0.01);
% Arrow at (-1,-1.5) pointing at 90 degrees (up)
\draw[-{Stealth[length=6pt]}] (-1,-1.6) -- ++(90:-0.01);
% Arrow at (0,-2) pointing at 180 degrees (left)
\draw[-{Stealth[length=6pt]}] (0.05,-2) -- ++(180:-0.01);
% Arrow at (1,-1.5) pointing at 270 degrees (down)
\draw[-{Stealth[length=6pt]}] (1,-1.4) -- ++(270:-0.01);
% Arrow at (0.5,-0.5) pointing at 315 degrees (down-right)
\draw[-{Stealth[length=6pt]}] (0.4,-0.4) -- ++(315:-0.01);
% Arrow at (0,-2.35) pointing at 180 degrees (left)
\draw[-{Stealth[length=6pt]}] (0.05,-2.37) -- ++(180:-0.01);
% Arrow at (0,-1.35) pointing at 0 degrees (right)
\draw[-{Stealth[length=6pt]}] (-0.05,-1.37) -- ++(0:-0.01);
% Arrow at (-0.5,0) pointing at 0 degrees (right)
\draw[-{Stealth[length=6pt]}] (-1,0) -- ++(0:-0.01);
% Arrow at (0.5,0) pointing at 0 degrees (right)
\draw[-{Stealth[length=6pt]}] (1,0) -- ++(0:0.01);
% Arrow at (0,-1) pointing at 0 degrees (right)
\draw[-{Stealth[length=6pt]}] (0.05,-1) -- ++(0:0.01);

\end{tikzpicture}
    \caption*{Embedded disks $[1,2,3,4,5]$, $[10,-2,7,5]$, $[2,9,4,8]$, and $[6,4,-7,-1]$.}
    \label{fig:5+4+4+42}
\end{figure}
%\end{proof}

\subsection{Fake Surfaces with Disconnected 1-Skeleta}
\label{sec:2-conn-components}
Using the classification of contractible fake surfaces up to complexity 5 from \cite{fagan2024classification}, it is easy to check that all satisfy one of the conditions of Lemma~\ref{lem:crl}. Thus, we can always 3-deform each surface to another of lower complexity. However, the classification is for fake surfaces with connected 1-skeleta. When we apply the reduction to lower complexity, we may get a fake surface whose 1-skeleton has 2 connected components. We henceforth use $\mathcal{F}_C(s,t)$ to refer to contractible fake surfaces of complexity $t$ whose 1-skeleton has $s$ connected components.

As discussed in the proof of Lemma~\ref{lem:crl}, only the trivial embedded disk move to embedded disks with one or two vertices from Case 1 may give a fake surface in $\mathcal{F}_C(2,t)$. In this case, the number of vertices will decrease by at least 2. Thus, it remains to address the case of $\mathcal{F}_C(2,t)$ for $t\leq 3$.

In \cite{ikeda1971acyclic}, it is shown the fake surfaces in $\mathcal{F}_C(2,1)$ satisfy the stable Andrews--Curtis and the Zeeman conjectures, so it remains to deal with the fake surfaces in $\mathcal{F}_C(2,2)$ and $\mathcal{F}_C(2,3)$. 

By Lemma 2.9 in \cite{koda2020shadows}, each surface $F\in\mathcal{F}_C(2,3)$ can be constructed in one of the following two ways:
\subsection*{Case 1}
In this case, we construct a fake surface $F=F_1\# F_2$ where $F_1\in\mathcal{F}_C(1,2)$ and $F_2\in\mathcal{F}(1,1)$ such that $H_1(F_2)=0$ and $H_2(F_2)=\mathbb{Z}$. Here, $\#$ represents the connected sum of two fake surfaces along their regions. It is straightforward to check that there are three possibilities for $F_2$. Their attaching maps are as follows, attached onto the 1-skeleton in Figure~\ref{fig:f2}: 
\begin{align*}
\text{Fake Surface 1: }& [-2,1,1],[1,2],[2]
\\
\text{Fake Surface 2: }& [-2,-1,2,1],[2],[1]
\\
\text{Fake Surface 3: }& [-2,1,-2,-1],[2],[1]
\end{align*}
\begin{figure}
    \centering
\begin{tikzpicture}[node distance=2cm, scale=2, transform shape] 
  \node[circle, fill=black, minimum size=4pt, inner sep=0pt] (A) {};
  % Self loop on A
  % two edges labeled 1 and 2
  \draw[postaction={decorate, decoration={markings, mark=at position 0.5 with {\arrow[scale=1.5]{>}}}}] 
    (A) to [out=135, in=225, looseness=30] node[midway, left, font=\tiny] {$1$} (A);
  \draw[postaction={decorate, decoration={markings, mark=at position 0.5 with {\arrow[scale=1.5]{>}}}}] 
    (A) to [out=45, in=-45, looseness=30] node[midway, right, font=\tiny] {$2$} (A);
\end{tikzpicture}
    \caption{1-Skeleton for $F_2$}
    \label{fig:f2}
\end{figure}

If $F_2$ is the first fake surface above, through computation about homology groups, the connected sum between $F_1$ and $F_2$ can not be done along the unique embedded disk $[2]$ of $F_2$. Thus $F=F_1\#F_2$ contains at least one embedded disk with one vertex such that we can reduce the complexity. When $F_2$ is either of the other two surfaces, we can always reduce the complexity since both of the surfaces contain two embedded disks.

\subsection*{Case 2} 
In this case, $F=F_1\cup F_2$, where $F_1$ is obtained by removing one disk in the interior of a 2-cell of one surface in $\mathcal{F}_C(1,3)$, and $F_2$ is obtained by attaching the boundary of a disk to the center of a M\"obius band. $F_1$ and $F_2$ are attached along their boundaries. We can collapse $F_2$ as follows: the 1-skeleton of $F_2$ is a circle disjoint from the 1-skeleton of $F_1$. Based on the construction of $F_2$, $F_2$ contains an embedded disk with the circle as the boundary and the $T$-bundle of the disk is nontrivial, as shown by the left blue disk in Figure~\ref{fig:F22}. We apply the inverse of the nontrivial embedded disk move, called the loop move in \cite{matveev2007algorithmic}, to the blue disk. The blue disk becomes a disk with one vertex and trivial $T$-bundle. Finally, we do the trivial embedded disk move to the new blue disk and remove the connected component of 1-skeleton of $F$.

We can apply the same argument to the two ways to construct the surfaces $F\in\mathcal{F}_C(2,2)$, following the construction in \cite{koda2020shadows}, and use the same method to reduce. 

Thus, we can use Lemma~\ref{lem:crl} to reduce each $F\in\mathcal{F}_C(1,t)(1<t\leq 5)$. When a reduction takes us outside $\mathcal{F}_C(1,t)$, we can use the above method to make it return to $\mathcal{F}_C(1,t)$ or reduce to a point. 

This completes the proof of Theorem~\ref{theorem:ac} in all cases.
\begin{figure}
    \centering
\begin{tikzpicture}
\draw (0,0) circle(1);

\draw[->,line width=1pt] (2,0)--(3,0);

\draw (5,0) circle(1);
\draw (6,0)..controls(9,3) and (9,-3)..(6,0);

\filldraw[opacity=0.3,blue] (0,0) circle(1);
\filldraw[opacity=0.3,blue] (5,0) circle(1);
\filldraw[opacity=0.3,red] (6,0)..controls(9,3) and (9,-3)..(6,0);

\node at (0,0){nontrivial};
\node at (5,0){trivial};
\node at (7.1,0){nontrivial};

\end{tikzpicture}
    \caption{Applying Inverse of Nontrivial Embedded Disk Move}
    \label{fig:F22}
\end{figure}

\section{Algorithmic 3-Deformations of Fake Surfaces}
\label{sec:3-def-fs-algos}
In this section, we provide pseudocode for implementing the embedded disk moves in Algorithms~\ref{alg:triv-bundle-move} and \ref{alg:nontriv-bundle-move} and a search algorithm to explore 3-deformations of fake surfaces in Algorithm~\ref{alg:fs-search}. For more details, see \cite{Fagan2025StableAC}. 

Our goal is to apply the embedded disk moves in order to execute a computer search to 3-deform fake surfaces to spines or ones of lower complexity.  While computer exploration has not yielded any new results so far, we plan to continue this approach in the future.

\begin{algorithm}
\caption{Trivial Embedded Disk Move}
\label{alg:triv-bundle-move}
\begin{algorithmic}[1]
\Function{TrivialEDMove}{$S, D, e$}
  \State Assert that $D$ is a trivial bundle
  \State Let $e$ be the front edge; list remaining edges as back edges
  \State Determine all perpendicular edge pairs using adjacency
  \State Introduce new vertical and horizontal edges
  \State Construct a replacement dictionary describing local edge substitutions
  \State Apply all replacements throughout the surface $S$
  \State Remove the original disk $D$
  \State Insert the new disks created by the move geometry
  \State Update edge indices and metadata
  \State \Return modified surface $S$
\EndFunction
\end{algorithmic}
\end{algorithm}

\begin{algorithm}
\caption{Nontrivial Embedded Disk Move}
\label{alg:nontriv-bundle-move}
\begin{algorithmic}[1]
\Function{NontrivialEDMove}{$S, D, e$}
  \State Assert that $D$ has a nontrivial $T$-bundle
  \State Let $e$ be the front edge; list back edges
  \State Compute perpendicular edge pairs, alternating orientation at each step
  \State Introduce new vertical, horizontal, and diagonal edges
  \State Build the replacement dictionary encoding the move
  \State Apply all replacements across $S$
  \State Remove the original disk $D$
  \State Insert the structured family of disks defined by the move
  \State Remove boundary edges and update metadata
  \State \Return modified surface $S$
\EndFunction
\end{algorithmic}
\end{algorithm}

\begin{algorithm}
\caption{Best-First Search Over Fake Surfaces}
\label{alg:fs-search}
\begin{algorithmic}[1]
\Function{Search}{$S_0, T, k$}
  \Comment{$T =$ max iterations, $k =$ moves sampled per expansion}
  \State Initialize priority queue $Q$
  \State Compute score of $S_0$ (complexity, number of nontrivial bundles)
  \State Insert $(\text{score}, S_0, \emptyset)$ into $Q$
  \State Track $S_{\mathrm{best}} \gets S_0$ and its score

  \For{$t = 1$ to $T$}
    \If{$Q$ is empty} \State \Return $S_{\mathrm{best}}$ \EndIf

    \State Pop best surface $(\_, S, H)$ from $Q$

    \State Determine all possible moves on $S$
    \State Sample up to $k$ moves

    \ForAll{sampled moves $m$}
       \State $S' \gets$ apply trivial or nontrivial embedded disk move $m$ to a copy of $S$
       \If{$S'$ is valid}
           \State Compute score of $S'$
           \State Insert $(\text{score}, S', H+[m])$ into $Q$
           \If{score($S'$) < score($S_{\mathrm{best}}$)}
               \State $S_{\mathrm{best}} \gets S'$
           \EndIf
       \EndIf
    \EndFor
  \EndFor

  \State \Return $S_{\mathrm{best}}$
\EndFunction
\end{algorithmic}
\end{algorithm}

\section*{Acknowledgments}
Y.Q. is supported by National Key R\&D Program of China (2024YFA1013202) and Nankai Zhide Foundation. Z.W. is partially supported by ARO MURI contract W911NF-20-1-0082.

\bibliographystyle{abbrv}
\bibliography{references}

\end{document}